\documentclass[a4paper, 11pt]{amsart}
\pdfoutput=1

\usepackage{amssymb}
\usepackage{amsfonts}
\usepackage{amsthm}
\usepackage{stmaryrd}
\usepackage{kpfonts} 
\usepackage[mathscr]{eucal}

\usepackage[english]{babel}

\usepackage[T1]{fontenc}
\usepackage[utf8]{inputenc}

\usepackage[margin=1in]{geometry}

\usepackage{tikz-cd}
\usepackage{tikz}
\usetikzlibrary[trees]
\usepackage{enumitem}
\usepackage{hyperref}
\usepackage{aliascnt}
\usepackage{adjustbox}
\usepackage{varioref}
\usepackage{xstring}
\usepackage{xcolor}
\usepackage{bbm}

\newtheorem{thmintro}{Theorem}

\newcommand{\inewtheorem}[2]{
	\newaliascnt{#1}{thmintro}
	\newtheorem{#1}[#1]{#2}
	\aliascntresetthe{#1}
}
\inewtheorem{propintro}{Proposition}
\inewtheorem{corintro}{Corollary}

\newtheorem{theorem}{Theorem}
\newcommand{\jnewtheorem}[2]{
	\newaliascnt{#1}{theorem}
	\newtheorem{#1}[#1]{#2}
	\aliascntresetthe{#1}
}
\numberwithin{theorem}{section}
\jnewtheorem{lemma}{Lemma}
\jnewtheorem{proposition}{Proposition}
\jnewtheorem{corollary}{Corollary}
\jnewtheorem{conjecture}{Conjecture}

\theoremstyle{definition}
\jnewtheorem{definition}{Definition}
\jnewtheorem{example}{Example}
\jnewtheorem{remark}{Remark}
\jnewtheorem{question}{Question}
\jnewtheorem{fact}{Fact}
\jnewtheorem{construction}{Construction}
\jnewtheorem{notation}{Notation}
\jnewtheorem{convention}{Convention}

\labelformat{theorem}{Theorem #1}
\labelformat{lemma}{Lemma #1}
\labelformat{proposition}{Proposition #1}
\labelformat{corollary}{Corollary #1}
\labelformat{definition}{Definition #1}
\labelformat{example}{Example #1}
\labelformat{section}{Section #1}
\labelformat{subsection}{Subsection #1}
\labelformat{equation}{(#1)}
\labelformat{remark}{Remark #1}
\labelformat{exercise}{Exercise #1}
\labelformat{question}{Question #1}
\labelformat{chapter}{Chapter #1}
\labelformat{construction}{Construction #1}
\labelformat{conjecture}{Conjecture #1}
\labelformat{notation}{Notation #1}
\labelformat{convention}{Convention #1}

\labelformat{thmintro}{Theorem #1}
\labelformat{corintro}{Corollary #1}
\labelformat{propintro}{Corollary #1}

\numberwithin{table}{subsection}
\labelformat{table}{Table #1}

\hypersetup{
	pdftitle={},
	pdfauthor={Gregoire Marc},
	bookmarksnumbered=true,     
	bookmarksopen=true,         
	bookmarksopenlevel=1,       
	colorlinks=true,   
	linkcolor = {cyan!60!black},   
	urlcolor = {cyan!60!black},   
	citecolor = {cyan!60!black},   
	pdfencoding=unicode,      
	pdfstartview=Fit,           
	pdfpagemode=UseOutlines,    
	pdfpagelayout=OneColumn,
	breaklinks=true,
	linktocpage
}

\def\DD{{\mathbb{D}}}
\def\EE{{\mathbb{E}}}
\def\FF{{\mathbf{F}}}

\def\NN{{\mathbb{N}}}

\def\sO{{\mathscr{O}}}
\def\sT{{\mathscr{T}}}

\def\sS{{\mathscr{S}}}

\def\cA{{\mathcal{A}}}
\def\cB{{\mathcal{B}}}
\def\cC{{\mathcal{C}}}
\def\cD{{\mathcal{D}}}
\def\cE{{\mathcal{E}}}

\def\cO{{\mathcal{O}}}

\def\cP{{\mathcal{P}}}
\def\cQ{{\mathcal{Q}}}
\def\cR{{\mathcal{R}}}
\def\cS{{\mathcal{S}}}
\def\cT{{\mathcal{T}}}
\def\cV{{\mathcal{V}}}
\def\cU{{\mathcal{U}}}

\newcommand\Seg[1]{\operatorname{Seg}_{#1}}
\newcommand\colim{\operatorname{colim}}
\newcommand\Cat{\operatorname{Cat}}
\newcommand\CAT{\widehat{\operatorname{Cat}}}
\newcommand\sCat{\operatorname{sCat}}

\newcommand\Id[1]{\operatorname{id}_{#1}}
\newcommand\Set{\operatorname{Set}}
\newcommand\opo[1]{{#1}^{\operatorname{op}}}
\newcommand\und[1]{\underline{#1}}
\newcommand\wspan{\operatorname{Span}_{\operatorname{L}}}
\newcommand\tspan{\operatorname{Span}_{\operatorname{L}}^{\operatorname{t}}}
\newcommand\fun{\operatorname{Fun}}

\newcommand\SSS{\mathbf{T}}

\newcommand\Del[1]{\Delta_{#1}}

\newcommand\Psh{\operatorname{Psh}}
\newcommand\Arr{\operatorname{Arr}}

\newcommand\num[1]{\underline{#1}}
\newcommand\sSet{\operatorname{sSet}}
\newcommand\Algg{\operatorname{Alg}}
\newcommand\lax{\operatorname{Lax}}
\newcommand\op{\operatorname{op}}
\newcommand\coll[1]{\operatorname{Coll}_{#1}}

\newcommand\operc[1]{\operatorname{Op}^{\operatorname{col}}_{#1}}
\newcommand\soperc[1]{{\operatorname{sOp}_{#1}^{\operatorname{col}}}}

\newcommand\oper[1]{\operatorname{Op}_{#1}}
\newcommand\soper[1]{\operatorname{sOp}_{#1}}

\newcommand\Lan{\operatorname{Lan}}
\newcommand\double{\operatorname{DoubleCat}}
\newcommand\suborb[1]{\operatorname{SubOrb}_{#1}}
\newcommand\MM{\mathbb{M}}
\newcommand\ccc[1]{\operatorname{C}_{#1}}

\newcommand\Adj{\operatorname{Adj}}

\newcommand\multi[1]{#1\text{-}\operatorname{Multicat}}

\newcommand\CSeg[1]{\operatorname{CSeg}_{#1}}
\newcommand\el{\operatorname{el}}
\newcommand\Span{\operatorname{Span}}
\newcommand\Fin{\operatorname{Fin}}

\title{Homotopy theory of simplicial parametrized operads}
\author{Gr\'{e}goire Marc}

\begin{document}
\maketitle
\begin{quote} 
\begin{center}
\textbf{Abstract} 
\end{center}
We define a generalization of (coloured) operads based on double lax functors and we construct a model structure on the associated category of generalized simplicial (coloured) operads. In particular, we obtain a model structure on the category of simplicial (coloured) $\sO$-operads of Nardin and Shah.  
\end{quote}
\tableofcontents

\section{Introduction}
Operads were introduced by May \cite{May} to study infinite loop spaces, and they play an important role in encoding algebraic structures in different areas of homotopy theory. In this paper, we are interested in coloured operads, which can be used to describe structures such as enriched categories or a ring together with a left module. Roughly speaking, a simplicial coloured operad $\cO$ is the data of a set of colours along with a simplicial set $\cO((x_1,\ldots,x_n),x_0)$ of ``multimorphisms from $(x_1,\ldots,x_n)$ to $x_0$'' for every list of colours $(x_0,\ldots,x_n)$. These simplicial sets are equipped with actions of the corresponding symmetric groups that ``permute the entries'' along with identities and an associative composition. The homotopy theory of simplicial operads has been studied in the non-coloured case in the work of Berger and Moerdijk \cite{BM03}, and a model structure on simplicial coloured operads has been constructed later on independently by Cisinski and Moerdijk \cite{cismo} and by Robertson \cite{Rob}. Simplicial coloured operads provide one model for the theory of \emph{$\infty$-operads} \cite{HA} that is known to be equivalent to the other models, combining the work of Cisinski and Moerdijk \cite{CisMoer,CisMoer2,cismo}, Barwick \cite{Bar2} and Chu, Haugseng and Heuts \cite{CHH}. 
 
In recent years, equivariant homotopy theory has become an increasingly active area of research, spurred in part by Hill, Hopkins, and Ravenel's resolution of the Kervaire invariant one problem \cite{HHR}. In particular, the theory of equivariant operads has gained particular attention through the development of \emph{$N_{\infty}$-operads} by Blumberg and Hill, which generalize $E_{\infty}$-operads to the equivariant setting. In contrast with the non-equivariant setting, where all $E_{\infty}$-operads are equivalent, 
$N_{\infty}$-operads are not all equivalent. Instead, they form a poset that can be described in terms of \emph{$G$-indexing systems}, following work of Blumberg and Hill 
\cite{BH1}, Bonventre and Pereira \cite{BonPer1}, and Rubin \cite{Rub1}. The study of the combinatorics of these $G$-indexing systems has recently become an active research area, for instance with the work of Balchin, Barnes, and Roitzheim \cite{BBR} or Balchin, MacBrough, and Ormsby \cite{BMO}. The study of equivariant operads is not limited to $N_{\infty}$-operads, see for example Hill's work on equivariant little disks operads \cite{Hill2022}. 

In \cite{BH1}, Blumberg and Hill defined the notion of \emph{graph equivalences of $G$-operads}, which is part of a model structure on the category of simplicial $G$-operads, as proven by Bonventre and Pereira \cite{BonPer1}. This graph model structure has then been extended to \emph{$G$-simplicial coloured operads} by Bonventre and Pereira in \cite{BonPer2}. Moreover, other models for equivariant operads, such as \emph{equivariant dendroidal sets} and \emph{equivariant dendroidal spaces}, have been developed and compared by Bonventre and Pereira in \cite{Per1,BonPer4,BonPer3}.

In the language of $\infty$-categories, based on Elmendorf's Theorem \cite{Elm}, equivariant homotopy theory can be formalized using the theory of \emph{parametrized $\infty$-categories}, as developed by Barwick, Dotto, Glasman, Nardin, and Shah in \cite{para, dede, Sha1}. This theory provides proper categorical foundations for higher equivariant algebra, and it has known recent further developments, for instance with the work of Cnossen, Lenz, and Linskens \cite{CLL1,CLL3,CLL2} or with the work of Hilman \cite{Hilman2024}. In particular, Lurie's theory of $\infty$-operads \cite{HA} admits a generalization in the equivariant world through the theory of \emph{parametrized $\infty$-operads} of Nardin and Shah \cite{DJ}, see also work of Barkan, Haugseng and Steinebrunner \cite{BHS}. In analogy with classical higher algebra, the theory of parametrized $\infty$-operads is expected to play a central role in \emph{higher equivariant algebra}. Progress in this direction has been made recently, for instance with the work of Stewart \cite{Nat4,Nat2,Nat3}.

The theory of equivariant operads has thus two facets: the classical theory of equivariant simplicial (or topological) operads and the theory of parametrized $\infty$-operads. Both these approaches have been studied with different perspectives, and it is natural to expect them to be the two sides of the same coin. One way to approach a comparison between these two theories is to introduce an intermediary player that shares one feature with each of these theories and that can be compared with both of them. A good candidate is the theory of \emph{simplicial parametrized operads} as introduced by Nardin and Shah in \cite{DJ}. The theory of simplicial parametrized operads and the theory of equivariant simplicial operads have one feature in common: they are both defined as model categories and can hopefully be compared via a Quillen adjunction. On the other hand, the theory of parametrized simplicial operads shares with parametrized $\infty$-operads that both these theories live in the parametrized world. 

\subsection{Main results}

In this paper, motivated by our overarching goal to provide a comparison of all models for equivariant operads, we construct a model structure on the category of simplicial (coloured) parametrized  operads. In order to achieve this, we consider a generalization of the theory of (coloured) symmetric operads that subsumes enriched categories, non-symmetric operads, and Nardin and Shah's parametrized operads. Our theory is based on the notion of \emph{$\FF$-multicategories} (see \ref{Fmulti} for a definition), that is a specialization of Leinster's generalized multicategories \cite{TL}. If $\cV$ is a cartesian category, then every $\FF$-multicategory $\MM$ gives rise to a notion of \emph{(coloured) $\MM$-operads in $\cV$} for which $\MM$ plays the role played by the symmetric groupoid $\Sigma_*$ for symmetric operads. Every $\MM$-operad in $\cV$ has an underlying \emph{$\MM$-collection in $\cV$}, and we prove, under some assumptions on $\MM$ and $\cV$, the existence of a composition product $\circ_{\MM}$ on the category $\coll{\MM}(\cV)$ of $\MM$-collections in $\cV$ for which algebras are precisely $\MM$-operads in $\cV$. More precisely, we prove the following theorem.

\begin{thmintro}(\ref{corocompprod})
For every $\FF$-multicategory $\MM$ that is target right fibrant in the sense of \ref{deftrf} and for every cartesian closed category $\cV$, the category $\coll{\MM}(\cV)$ of $\MM$-collections in $\cV$ is endowed with the structure of a monoidal category for which the category of algebras is precisely the category of $\MM$-operads in $\cV$.
\end{thmintro}

Using this generalization of the composition product, and the fact that $\coll{\MM}(\sSet)$ can be endowed with a projective model structure, we prove the following theorem. 

\begin{thmintro}(\ref{THmodop})\label{THwcol}
For every target right fibrant $\FF$-multicategory $\MM$, the category $\soper{\MM}$ admits a model structure that is right lifted from $\coll{\MM}(\sSet)$ with respect to the free-forgetful adjunction
$$
F_{\MM} \colon \coll{\MM}(\sSet) \rightleftarrows \soper{\MM} \colon U_{\MM}.
$$
\end{thmintro}

One main feature of our notion of generalized operads is that it provides a convenient way to encode coloured operads. For every \emph{$\MM$-set of colours} $A$ (in the sense of \ref{defcolour}), $A$-coloured operads are by definition $\MM^{A}$-operads for some $\FF$-multicategory $\MM^{A}$. Using this description of coloured $\MM$-operads, we prove the following theorem.

\begin{thmintro}(\ref{mainTH})
For every $\FF$-multicategory $\MM$ that is target right fibrant, the category of simplicial coloured $\MM$-operads admits a model structure for which
the weak equivalences are the fully faithful essentially surjective morphisms and for which the fibrations are the isofibrations (See \ref{weakfibop} for all the terminology used here).
\end{thmintro}

We can specialize this theorem to obtain a model structure on Nardin and Shah's simplicial parametrized operads.

\begin{corintro}(\ref{modparacop})
For every category $\sO$ and every orbital subcategory $\sT$ of $\sO$, the category $\soperc{(\sO,\sT)}$ of simplicial coloured $(\sO,\sT)$-parametrized operads admits a model structure for which
the weak equivalences are the fully faithful essentially surjective morphisms and for which the fibrations are the isofibrations (respectively defined in \ref{weakfibop}).
\end{corintro}

\subsection{Structure of the paper}

In \ref{sec2}, we define the notion of $\FF$-multicategories as well as our notion of generalized operads. In \ref{sec3},  we construct, for a target right fibrant $\FF$-multicategory $\MM$, a convolution product on the category of $\MM$-collections for which the algebras are precisely the $\MM$-operads. In \ref{sec4}, for a target right fibrant $\FF$-multicategory $\MM$, we construct a model structure on the category of simplicial (coloured) $\MM$-operads. Finally, in \ref{sec5}, we define a notion of \emph{complete $\MM$-Segal spaces}, which we conjecture to be a model for simplicial $\MM$-operads.

\subsection{Related work}

For a finite group $G$ and a $G$-transfer system $\cT$, the category $\soper{(\sO_G,\sO_{\cT})}$ defined in \ref{exgenop} should be equivalent to the category of \emph{simplicial $\cT$-partial genuine equivariant operads} of Bonventre and Pereira \cite{BonPer1}, however, such a comparison has not been made yet. In particular, we expect that the model structure on $\soper{(\sO_G,\sO_{\cT})}$ that we obtain in \ref{exgenmod} corresponds to their model structure. 

This paper is based on the idea that the theory of double categories can be used to implement operadic theories. This idea is not new, and such considerations have already been made by Pisani in \cite{Pis} and by Haugseng in \cite{Haug1}.

\subsection{Notations and conventions}

We fix three Grothendieck universes $\mathbf{U}\in \mathbf{V} \in \mathbf{W}$ that contain the first infinite ordinal. A set is \emph{small} if it is contained in $\mathbf{U}$ and \emph{large} if it is contained in $\mathbf{W}$. We will also use an analogous naming convention for categories. In what follows, we denote by $\Set$ the category of small sets, by $\sSet$ the category of small simplicial sets, by $\Cat$ the category of small categories, and by $\CAT$ the large category of categories.

\subsection{Acknowledgements}
I would like to thank my PhD supervisor, Magdalena Kędziorek, for her continuous support and for her valuable comments on earlier drafts of this paper. The author would like to thank the Isaac Newton Institute for Mathematical Sciences, Cambridge, for support and hospitality during the programme Equivariant homotopy theory in context, where work on this paper was undertaken. This work was supported by EPSRC grant EP/Z000580/1. During the writing of this paper, the author was supported by an NWO grant Vidi.203.004.

\section{\texorpdfstring{$\MM$}{M}-operads}\label{sec2}

\subsection{\texorpdfstring{$\FF$}{F}-multicategories}\label{Fmulti}
The goal of this section is to define the notion of $\FF$-multicategories as a specialization of Leinster's generalized multicategories \cite{TL}.
\begin{definition}[\protect{\cite[Section 3]{Ben1}}]
Let $\cC$ be a category. The \emph{completion by finite coproducts of $\cC$} is the category $\FF \cC$ whose objects are sequences $\underline{X}=(X_1,\ldots, X_n)$ of objects of $\cC$ and, given two such sequences $\underline{X}=(X_1,\ldots, X_n)$ and $\underline{Y}=(Y_1,\ldots, Y_m)$, whose set of morphisms $\FF \cC(\underline{X},\underline{Y})$ is 
\[
\coprod_{\alpha\colon \num{n}\to \num{m}}\prod_{i=1}^{n} \cC(X_{i}, Y_{\alpha(i)}),
\]
where $\num{n}$ denotes the finite set $\{1,\ldots ,n\}$.
The composition of morphisms in $\FF \cC$ is naturally induced by the composition of maps of finite sets and by the composition of morphisms in $\cC$. 
\end{definition}

\begin{proposition}\label{univfin}
For every category $\cC$, the category $\FF \cC$ admits finite coproducts and satisfies the following universal property: for every category $\cD$ with finite coproducts, the canonical inclusion functor $\cC \to \FF \cC $ induces an equivalence of categories
$$
\fun_{\sqcup} (\FF \cC,\cD)\to \fun (\cC,\cD)
$$
where $\fun_{\sqcup} (\FF \cC,\cD)$ denotes the full subcategory of $\fun (\FF \cC ,\cD)$ spanned by these functors that preserve finite coproducts.
\begin{proof}
This proposition is classical. The coproducts in the category $\FF \cC$ are given by the concatenation of sequences of objects in $\cC$.
\end{proof}
\end{proposition}

\begin{remark}\label{opof}
For every category $\cC$, the category $\opo{(\FF\cC)}$ is the completion by finite products of $\opo{\cC}$. In particular, for every functor $F\colon \opo{\cC} \to \cV$ where $\cV$ is a category with finite products, $F$ can be extended into a functor $\widetilde{F}\colon \opo{(\FF \cC)} \to \cV$ that preserves finite products defined explicitly on objects by the formula
$$
\widetilde{F}(X_1,\ldots,X_n)=\prod_{i=1}^n F(X_i).
$$
\end{remark}

\begin{proposition}
The functor $\FF \colon \Cat \to \Cat$ is naturally endowed with the structure of a cartesian monad.
\begin{proof}
For every category $\cC$, we obtain a functor $\mu_\FF \colon \FF^2 \cC  \to \FF\cC$ induced by the identity $\Id{\cC}\colon \FF \cC \to \FF \cC$ using \ref{univfin}. A direct verification shows that the natural transformation $\mu \colon \FF^2 \cC \to \FF \cC$ thus obtained, along with the inclusion functor $\cC \hookrightarrow \FF \cC$, induce the structure of a cartesian monad on $\FF$. 
\end{proof}
\end{proposition}

\begin{remark}
Note that the structure of a monad on $\FF \cC$ exists only because we chose a strict model for $\FF \cC$. For a different model, we would only obtain the structure of a pseudo monad on $\FF$. 
\end{remark}

Specializing \cite[Definition 4.2.2]{TL} at the monad $\FF$, we obtain a notion of $\FF$-multicategories as well as a notion of multifunctors between $\FF$-multicategories. We denote by $\multi{\FF}$ the category of $\FF$-multicategories and multifunctors. Recall that an $\FF$-multicategory $\MM$ is given by the following data:
\begin{enumerate}
    \item a small category $\MM_1$ of \emph{multimorphisms} and a small category $\MM_0$ of \emph{objects};
    \item a \emph{source functor} $s\colon \MM_1 \to \FF \MM_0$ and a \emph{target functor} $t\colon \MM_1 \to \MM_0$;
    \item an \emph{identity functor} $1_{\MM}\colon \MM_0 \to \MM_1$ and a composition functor $\odot\colon \MM_1 \times_{\FF \MM_0} \FF \MM_1 \to \MM_1$.
\end{enumerate}
These data are required to satisfy the following compatibilities:
\begin{enumerate}
    \item the composition of $1_{\MM}$ and $s$ is the canonical inclusion $\MM_0 \hookrightarrow \FF \MM_0$ and $1_{\MM}$ is a section of $t$;
    \item the composition of the composition functor $\odot\colon \MM_1 \times_{\FF \MM_0} \FF \MM_1 \to \MM_1$ with the target functor $t\colon \MM_1 \to \MM_0$ is the composition
$$
\MM_1 \times_{\FF \MM_0} \FF \MM_1 \overset{p_1}{\longrightarrow} \MM_1 \overset{t}{\longrightarrow} \MM_0
$$
and the composition of the functor $\odot\colon \MM_1 \times_{\FF \MM_0} \FF \MM_1 \to \MM_1$ with the source functor $s\colon \MM_1 \to \MM_0$ is the composition
$$
\MM_1\times_{\FF \MM_0} \FF \MM_1 \overset{p_2}{\longrightarrow} \FF \MM_1 \overset{\FF s}{\longrightarrow} \FF^2 \MM_0 \overset{\mu_{\FF}}{\longrightarrow} \FF \MM_0;
$$
    \item (associativity) the following diagram commutes
$$
\begin{tikzcd}
\MM_1 \times_{\FF \MM_0} \FF \MM_1 \times_{\FF^2 \MM_0} \FF^2 \MM_1 \arrow[d,"\odot\times_{\Id{\FF^2\MM_1}}\Id{\FF^2\MM_1}"'] \arrow[r,"\simeq"] & \MM_1 \times_{\FF \MM_0}\FF\left(\MM_1 \times_{\FF \MM_0} \FF \MM_1 \right) \arrow[d,"\Id{\MM_1}\times_{\Id{\FF \MM_0}}\FF \odot"]  \\
\MM_1 \times_{\FF^2 \MM_0}\FF^2 \MM_1 \arrow[d,"\simeq"'] & \arrow[d,"\odot"] \MM_1 \times_{\FF \MM_0}\FF \MM_1 \\
\MM_1 \times_{\FF \MM_0} \FF \MM_1 \arrow[r,"\odot"'] & \MM_1;
\end{tikzcd}
$$
    \item (unitality) the following diagrams commute
$$
\begin{tikzcd}[column sep=huge]
\MM_1\times_{\FF \MM_0} \FF\MM_0 \arrow[r,"\Id{\MM_1}\times_{\Id{\MM_0}}\FF 1_{\MM}"] \arrow[dr,"\simeq"'] & \arrow[d,"\odot"] \MM_1 \times_{\FF \MM_0} \FF \MM_1 &  \MM_0 \times_{\MM_0} \MM_1 \arrow[r,"1_{\MM}\times_{} "] \arrow[rd,"\simeq"'] & \arrow[d,"\odot"]  \MM_1 \times_{\FF \MM_0} \FF \MM_1 \\
& \MM_1 & &  \MM_1.
\end{tikzcd}
$$
\end{enumerate}

Recall now from \cite{Ehr1} the notion of a double category.
\begin{definition}\label{defdoub}
A \emph{double category} is a category internal to $\CAT$. Equivalently, a double category is an $\Id{}$-multicategory in the sense of \cite[Definition 4.2.2]{TL} where $\Id{}\colon \CAT \to \CAT$ is the identity monad. A double category is \emph{small} if it is a category internal to $\Cat$. We denote by $\double$ the category of small double categories.
\end{definition}

\begin{remark}
Note that \ref{defdoub} is the definition of strict double categories. For the rest of this document, we will only consider strict double categories, and we will always call them double categories.
\end{remark}

\begin{proposition}
The functor $\FF \colon \Cat \to \Cat$ is naturally endowed with the structure of a cartesian oplax morphism of monad $\FF \colon (\Cat,\FF) \to (\Cat, \Id{})$ in the sense of \cite[Section 6.7]{TL}.
\begin{proof}
The functor $\FF$ is cartesian and the structure of colax functor on $\FF$ is provided by the multiplication of $\FF$.
\end{proof}
\end{proposition}

Using \cite[Section 6.7]{TL} on the colax functor $\FF \colon (\Cat,\FF) \to (\Cat, \Id{})$, we obtain a functor 
$$
\FF_* \colon \multi{\FF} \to \double
$$
that sends an $\FF$-multicategory $\MM$ to the small double category $\FF \MM$ defined by $(\FF \MM)_1=\FF \MM_1$ and $(\FF \MM)_0=\FF \MM_0$ and where the source and target functors are respectively given by the composition $\FF \MM_1 \overset{\FF s_{\MM}}{\longrightarrow} \FF^2\MM_0 \overset{\mu_{\FF}}{\longrightarrow}\FF \MM_0$ and by $\FF t_{\MM}\colon \FF \MM_1 \to \FF \MM_0$, where the horizontal composition functor is given by $\FF \odot \colon \FF \MM_1 \times_{\FF \MM_0} \FF \MM_1 \to \FF \MM_1$ using that $\FF\left(\FF\MM_1 \times_{\FF \MM_0} \MM_1 \right)\simeq \FF \MM_1 \times_{\FF \MM_0} \FF \MM_1$ and where the identity functor is $\FF 1_{\MM}\colon \FF \MM_0 \to \FF \MM_1$.

\subsection{\texorpdfstring{$\MM$}{M}-operads}
In this section, we define our generalization of operads based on the notion of $\FF$-multicategories. 

\begin{remark}\label{monasdoub}
Recall that every monoidal category $(\cC,\otimes)$ can be seen as a weak double category with a trivial category of vertical morphisms and whose category of horizontal morphisms and squares is $\cC$ itself. The horizontal composition is given by the monoidal structure of $(\cC,\otimes)$.
\end{remark}

We now recall from \cite{GrPa} the definition of lax functors between weak double categories. We only recall the definition when the target is a monoidal category and when the source is a (strict) double category.
\begin{definition}\label{deflax}
Let $\DD$ be a double category and let $(\cC,\otimes_{\cC})$ be a monoidal category. A \emph{lax functor $F\colon \DD \to \cC$} is a functor $F\colon \DD_1 \to \cC$ along with a natural transformation $\alpha\colon \otimes_{\cC}\circ (F\times F) \Rightarrow F \circ \odot$ of functors $\DD_1 \times_{\DD_0} \DD_1 \to \cC$ and a natural transformation $\eta \colon 1_{\cC} \Rightarrow F \circ 1_{\DD}$ of functors $\DD_0 \to \cC$ satisfying the following compatibilities:
\begin{enumerate}
    \item for every $(f_1,f_2,f_3)$ in $\DD_1 \times_{\DD_0} \DD_1 \times_{\DD_0} \DD_1$, the following diagram commutes
$$
\begin{tikzcd}[column sep=large]
F(f_1) \otimes_{\cC} F(f_2) \otimes_{\cC} F(f_3) \arrow[r,"\alpha\otimes_{\cC} \Id{F(f_3)}"] \arrow[d,"\Id{F(f_1)}\otimes_{\cC} \alpha"'] & F(f_1\odot f_2)\otimes_{\cC} F(f_3) \arrow[d,"\alpha"] \\
F(f_1) \otimes_{\cC} F(f_2 \odot f_3) \arrow[r,"\alpha"'] & F(f_1 \odot f_2 \odot f_3);
\end{tikzcd}
$$
    \item for every $X$ in $\DD_0$ and every $f$ in $\DD_1$, the following triangles commute
$$
\begin{tikzcd}
F(f) \arrow[r,"\Id{F(f)}\otimes_{\cC} \eta"] \arrow[rd,"="'] & F(f)\otimes_{\cC} F(1_X) \arrow[d,"\alpha"] & & F(f) \arrow[r,"\eta \otimes_{\cC} \Id{F(f)}"] \arrow[rd,"="'] & F(1_X)\otimes_{\cC} F(f) \arrow[d,"\alpha"] \\
& F(f\odot 1_X) & & & F(1_X \odot f).
\end{tikzcd}
$$
\end{enumerate}
A \emph{vertical natural transformation $\theta\colon F \Rightarrow G$} between two lax functors $F,G\colon \DD \to \cC$ is a natural transformation $\theta\colon F \Rightarrow G$ between the underlying functors $F,G\colon \DD_1 \to \cC$ satisfying the following conditions:
\begin{enumerate}
    \item for every $(f,g)$ in $\DD_1 \times_{\DD_0} \DD_1$, the following square commutes
$$
\begin{tikzcd}
F(f) \otimes_{\cC} F(g) \arrow[r,"\theta_f \otimes \theta_g "] \arrow[d,"\alpha"'] & G(f) \otimes_{\cC} G(g) \arrow[d,"\alpha"] \\
F(f\odot g) \arrow[r,"\theta_{f\odot g}"'] & G(f\odot g);
\end{tikzcd}
$$
    \item for every $X$ in $\DD_0$, the following triangle commutes
$$
\begin{tikzcd}[row sep=small]
& F(1_X) \arrow[dd,"\theta_{1_X}"] \\
1_{\cC}\arrow[ru,"\eta"] \arrow[rd,"\eta"'] & \\
& G(1_X).
\end{tikzcd}
$$
\end{enumerate}
We denote by $\lax (\DD,\cC) $ the category of lax functors $F\colon \DD \to \cC$ and vertical natural transformations thus obtained. 
\end{definition}

\begin{definition}
Let $\MM$ be an $\FF$-multicategory and let $\cV$ be a category with finite products.
An \emph{$\MM$-collection in $\cV$} is a functor $S \colon \opo{\MM_1} \to \cV$.
An \emph{$\MM$-operad in $\cV$} is a lax functor $\cO \colon \opo{\FF \MM} \to \cV$ in the sense of \ref{deflax} such that the underlying functor $\cO \colon \opo{\FF \MM_1} \to \cV$ preserves finite products. Explicitly, an $\MM$-operad $\cO$ in $\cV$ is the data of an $\MM$-collection $\cO\colon \opo{\MM_1} \to \cV$ in $\cV$ along with:
\begin{enumerate}
    \item identities given by a natural transformation $1_{X}\colon 1_{\cV} \to \cO(1_{X})$ of functors $\opo{\MM_0} \to \cV$ where the right hand side is the composition
$$
\opo{\MM_0} \overset{1_{\MM}}{\longrightarrow} \opo{\MM_1} \overset{\cO}{\longrightarrow} \cV,
$$
and where the left hand side is the functor constant on $1_{\cV}$;

    \item a composition law given by a natural transformation $\circ \colon \cO(f)\times \cO(f_1)\times \cdots \times \cO (f_n) \to \cO \left (f\odot(f_1,\ldots, f_n)\right)$ of functors 
$\opo{\left(\MM_1 \times_{\FF \MM_0} \FF \MM_1 \right)} \to \cV$ where the left hand side is the composition 
$$
\opo{\left(\MM_1 \times_{\FF \MM_0} \FF \MM_1 \right)} \overset{}{\longrightarrow} \opo{\left(\MM_1 \times \FF \MM_1\right)} \overset{\cO \times \widetilde{\cO}}{\longrightarrow} \cV \times \cV \overset{}{\longrightarrow} \cV, 
$$
with $\widetilde{\cO}$ given by \ref{opof}, and where the right hand side is the composition
$$
\opo{\left(\MM_1 \times_{\FF \MM_0} \FF \MM_1 \right)} \overset{\opo{\odot}}{\longrightarrow} \opo{\MM_1} \overset{\cO}{\longrightarrow} \cV.
$$
\end{enumerate}
These data are required to satisfy the following compatibilities:
\begin{enumerate}
    \item (associativity) for every
$$
(f,(f_1,\ldots,f_n),((f_1^1,\ldots,f_1^{m_1}),\ldots,(f_n^1,\ldots,f_n^{m_n}))) \in \MM_1 \times_{\FF \MM_0} \FF \MM_1 \times_{\FF^2 \MM_0} \FF^2 \MM_1,
$$
the following diagram commutes
$$
\begin{tikzcd}
\arrow[r,"\simeq"] \cO(f) \times \prod_{i=1}^n\left( \cO(f_i)\times \prod_{j=1}^{m_j} \cO(f_i^j) \right) \arrow[d,"\Id{}\times \prod_{i=1}^n\odot"'] &  \arrow[d,"\odot\times \Id{}"]  \cO(f) \times \left(\prod_{i=1}^n \cO(f_i) \right) \times \left(\prod_{i=1}^n \prod_{j=1}^{m_i}\cO(f_i^j) \right)   \\
 \cO(f)\times \prod_{i=1}^n \cO(f_i\odot(f_i^1,\ldots,f_i^{m_i})) \arrow[d,"\odot"']&   \cO(f\odot (f_1,\ldots, f_n)) \times \left(\prod_{i=1}^n \prod_{j=1}^{m_i}\cO(f_i^j) \right) \arrow[d,"\odot"]  \\
\cO(f\odot(f_1\odot(f_1,\ldots,f_1^{m_1}),\ldots,f_n\odot(f_n^1,\ldots,f_n^{m_n}))) \arrow[r,"="'] & \cO((f\odot (f_1,\ldots,f_n))\odot (f_1^1,\ldots,f_n^{m_n}))  ;
\end{tikzcd}
$$
    \item (unitality) for every $f$ in $\MM_1$, if we denote by $(x_1,\ldots,x_n)$ the source of $f$, then the following diagram commutes
$$
\begin{tikzcd}[column sep=large]
\cO(f) \arrow[r,"(\Id{}\text{,} 1_{x_1}\text{,}\ldots \text{,}1_{x_n})"] \arrow[rd,"="'] & \cO(f)\times \cO(1_{x_1}) \times \cdots \cO(1_{x_n}) \arrow[d,"\odot"] \\
 & \cO(f\odot(1_{x_1},\ldots,1_{x_n})),
\end{tikzcd}
$$
and if we denote by $x$ the target of $f$, then the following diagram commutes
$$
\begin{tikzcd}[column sep=large]
\cO(f) \arrow[r,"(\Id{}\text{,}1_x)"] \arrow[rd,"="'] & \cO(f)\times \cO(1_{x}) \arrow[d,"\odot"] \\
 & \cO(1_x \odot f).
\end{tikzcd}
$$
\end{enumerate}
A \emph{morphism of $\MM$-operads $\phi\colon \cO \to \cP$ in $\cV$}  is a vertical natural transformation $\phi \colon \cO \Rightarrow \cP$. Explicitly, a morphism $\phi\colon \cO \to \cP$ of $\MM$-operads in $\cV$ is a natural transformation $\phi\colon \cO \Rightarrow \cP$ of functors $\opo{\MM_1} \to \cV$ such that the following diagrams commute:
\begin{enumerate}
    \item 
$$
\begin{tikzcd}[column sep=huge]
\cO(f)\times \cO(f_1)\times \cdots \times \cO(f_n) \arrow[r,"\phi(f)\times \prod_{i=1}^n \phi(f_i)"] \arrow[d,"\circ_{\cO}"'] & \cP(f) \times \cP(f_1)\times \cdots \times \cP(f_n)    \arrow[d,"\circ_{\cP}"] \\
\cO(f\odot (f_1,\ldots,f_n)) \arrow[r,"\phi(f\odot(f_1\text{,}\ldots\text{,}f_n))"'] & \cP(f\odot (f_1,\ldots,f_n));
\end{tikzcd}
$$
    \item
$$
\begin{tikzcd}[row sep=small]
& \cO(1_X) \arrow[dd,"\phi(1_X)"] \\ 
1_{\cV} \arrow[ru,"1_X"] \arrow[rd,"1_X"'] & \\
& \cP(1_X).
\end{tikzcd}
$$
\end{enumerate}
In what follows, we denote by $\oper{\MM}(\cV)$ the category of $\MM$-operads in $\cV$ thus obtained. Note that, by definition, $\oper{\MM}(\cV)$ is a full subcategory of $\lax (\opo{\FF \MM},\cV)$.
\end{definition}

\begin{example}\label{exsymmop}
Consider the $\FF$-multicategory $\und{\Sigma_*}$ whose category of multimorphisms is the symmetric groupoid $\Sigma_*=\bigsqcup_{n\ge 0}\operatorname{B}\Sigma_n$, whose category of objects is the final category $*$, whose source and target functors are the inclusion $\Sigma_* \hookrightarrow \FF *$ and the unique functor $\Sigma_* \to *$ respectively, and whose composition functor $\odot \colon \Sigma_* \times_{\FF *}\FF \Sigma_* \to \Sigma_*$ sends an object $(n,(k_1,\ldots,k_n))$ to $k_1+ \cdots+ k_n$ and sends a morphism $(\delta,(\eta_1,\cdots,\eta_n))$ with $\delta$ in $\Sigma_n$ and $\eta_i$ in $\Sigma_{k_i}$ to the permutation 
$$
\delta_{k_1,\ldots,k_n}\circ \left( \bigsqcup_{i=1}^k \eta_i \right),
$$
with $\delta_{k_1,\ldots,k_n}$ the block permutation. The category of $\und{\Sigma_*}$-operads in $\cV$ is the category of symmetric operads in $\cV$.
\end{example}

\begin{example}\label{exnonsym}
Consider the $\FF$-multicategory $\und{\NN}$ whose category of multimorphisms is the set $\NN$ of non negative integers (seen as a discrete category), whose category of objects is $*$, whose source and target functors are the inclusion $\NN \hookrightarrow \FF *$ and the unique functor $\NN \to *$ respectively, and whose composition functor $\NN \times _{\FF * } \FF \NN$ sends an objects $(n,(k_1,\ldots,k_n))$ to $k_1+\cdots+k_n$.
The category of $\und{\NN}$-operads in $\cV$ is the category of non-symmetric operads in $\cV$.
\end{example}

\begin{remark}
Note that $\und{\NN}$ is actually a classical multicategory seen here as a discrete $\FF$-multicategory.
\end{remark}

\begin{proposition}\label{laxfunc}
For every multifunctor $\phi \colon \MM \to \NN$ between $\FF$-multicategories and for every functor $\psi \colon \cV \to \cU$ between categories with finite products that preserves finite products, we obtain by precomposition and post-composition functors $\phi^*\colon \oper{\NN}(\cV) \to \oper{\MM}(\cV)$ and $\oper{\NN}(\psi)\colon \oper{\NN}(\cV) \to \oper{\NN}(\cU)$. In particular, for every category with finite products $\cV$, we obtain a functor
$$
\oper{(-)}(\cV)\colon \opo{\multi{\FF}} \to \CAT.
$$
\begin{proof}
This follows directly from the description of $\MM$-operads in terms of lax functors.
\end{proof}
\end{proposition}

\subsection{Coloured \texorpdfstring{$\MM$}{M}-operads}
The goal of this section is to define the notion of coloured $\MM$-operads. Recall first the definition of the Grothendieck construction from \cite{SGA1}.

\begin{definition}\label{defgrotcons}
Consider a functor $F\colon \opo{\cC} \to \CAT$. The \emph{Grothendieck construction of $F$} is the category $\int_{\cC} F$ whose objects are pairs $(X,x)$ with $X$ in $\cC$ and $x$ in $F(X)$ and whose morphims $(X,x) \to (Y,x)$ are pairs $(f,\phi)$ with $f\colon X \to Y$ a morphism in $\cC$ and $\phi\colon x \to F(f)(y)$ a morphism in $F(X)$. We obtain moreover a functor $\pi_F\colon \int_{\cC} F \to \cC$ that sends $(X,x)$ to $X$. If $F$ is a presheaf of small sets on $\cC$, \emph{the category of elements of $F$} denoted by $\el(F)$ is the Grothendieck construction of the composition
$$
\opo{\cC} \overset{F}{\longrightarrow} \Set \hookrightarrow \CAT
$$
where the right hand side functor sends a set to the corresponding discrete category.
\end{definition}

\begin{remark}
Explicitly, the category of elements of a presheaf $F\colon \opo{\cC} \to \Set$ is the small category whose objects are pairs $(X,x)$ with $X$ in $\cC$ and $x$ in $F(X)$ and where a morphism $(X,x) \to (X,y)$ is a morphism $f\colon X \to Y$ in $\cC$ such that $F(f)(y)=x$.
\end{remark}

\begin{definition}\label{defcolour}
For every $\FF$-multicategory $\MM$, the category of \emph{$\MM$-sets of colours} is the category $\Psh (\MM_0)$ of presheaves of small sets on $\MM_0$. 
\end{definition}

\begin{construction}
For every $\FF$-multicategory $\MM$ and every functor $F:\cA \to \MM_0$ with $\cA$ a small category, we define a new $\FF$-multicategory $F^*\MM$ whose category of objects $(F^*\MM)_0$ is $\cA$ and whose category of multimorphisms $(F^*\MM)_1$ is $\MM_1 \times_{\left(\MM_0 \times \FF \MM_0\right)} \left( \cA \times\FF \cA \right)$, whose composition functor 
$$
\odot_{F^* \MM} \colon  F^*\MM_1\times_{\FF \left(F^*\MM_0 \right)}\FF\left(F^*\MM_1\right) \to F^* \MM_1
$$ 
is the functor
$$
\left( \MM_1 \times_{\FF \MM_0} \FF \MM_1 \right) \times_{\left( \MM_0 \times \FF \MM_1 \times \FF^2 \MM_0 \right)} ( \cA\times \FF \cA \times \FF^2 \cA) \overset{\circ \times \mu_{\FF}\times \Id{\cA} }{\longrightarrow} \MM_1 \times_{\left(  \MM_0 \times \FF \MM_0 \right)}\left( \cA \times \FF \cA \right)
$$
using the isomorphism $F^*\MM_1 \times_{\FF\left(F^*\MM_0 \right)}\FF \left(F^*\MM_1\right) \simeq \left( \MM_1 \times_{\FF \MM_0} \FF \MM_1 \right) \times_{\left( \MM_0 \times \FF \MM_1 \times \FF^2 \MM_0 \right)} ( \cA\times \FF \cA \times \FF^2 \cA)$, and whose identity functor $1\colon F^*\MM_0 \to F^*\MM_1$ is the functor $\cA \to \MM_1 \times_{\left( \MM_0 \times \FF \MM_0\right)} \left(  \cA \times \FF \cA \right)$ given component-wise by the composition $ \cA \overset{F}{\longrightarrow} \MM_0 \overset{1}{\longrightarrow} \MM_1$, by the inclusion $ \cA \hookrightarrow \FF \cA$ and by the identity $\Id{ \cA}$ respectively.
In particular, we obtain a functor $(-)^*\MM \colon \Cat_{/\MM_0} \to \multi{\FF}$ that sends a functor $F\colon \cA \to \MM_0$ to the $\FF$-multicategory $F^* \MM$ and that sends a commutative triangle of functors
$$
\begin{tikzcd}[column sep=small]
\cA \arrow[rr,"H"] \arrow [rd,"F"'] & & \cB \arrow[ld,"G"] \\
& \MM_0 &
\end{tikzcd}
$$ 
to the multifunctor $F^*\MM \to G^*\MM$ given on objects by $H$ and on morphisms by the functor $(H^* \MM)_1\colon (F^* \MM)_1 \to (G^*\MM)_1$ obtained by the following diagram
$$
\begin{tikzcd}
(F^*\MM)_1 \arrow[d,"(H^* \MM)_1"] \arrow[dd,bend right=50,""'] \arrow[r] & \cA \times \FF \cA \arrow[d," H \times \FF H"] \\
(G^*\MM)_1 \arrow[d] \arrow[r] & \cB \times \FF \cB \arrow[d," G \times \FF G"] \\
\MM_1 \arrow[r,"s\times t"'] &  \MM_0 \times \FF \MM_0.
\end{tikzcd}
$$  
In what follows, we denote by $\ccc{\MM}$ the composition $\Psh(\MM_0) \overset{\el}{\longrightarrow} \Cat_{/\MM_0} \overset{(-)^*\MM}{\longrightarrow} \multi{\FF}$ and we denote by $\MM^A$ the image of a $\MM$-set of colours $A$ under the functor $\ccc{\MM}$.
\end{construction}

\begin{definition}\label{colcol}
Let $A$ be an $\MM$-set of colours. An \emph{$A$-coloured $\MM$-collection in $\cV$} is an $\MM^A$-collection in $\cV$. 
Similarly, an \emph{$A$-coloured $\MM$-operad in $\cV$} is an $\MM^A$-operad in $\cV$. We denote respectively by $\coll{\MM}^A(\cV)$ and $\oper{\MM}^A(\cV)$ the categories of $A$-coloured collections and $A$-coloured operads in $\cV$.
Consider the functor $\Phi_{\MM}  \colon \opo{\Psh (\MM_0)} \to \CAT$ that sends $A$ to $\oper{\MM}^A(\cV)$, obtained as the composition
$$
\opo{\Psh(\MM_0)} \overset{\ccc{\MM}}{\longrightarrow} \opo{\left(\multi{\FF}\right)} \overset{\oper{(-)}(\cV)}{\longrightarrow} \CAT.
$$
We denote by $\operc{\MM}(\cV)$ the Grothendieck construction of $\Phi_{\MM}$. The category $\operc{\MM}(\cV)$ is the category of \emph{coloured $\MM$-operads in $\cV$}.
\end{definition}

\begin{example}
If we consider the $\FF$-multicategory $\und{\Sigma_*}$ of \ref{exsymmop}, then the category $\operc{\und{\Sigma_*}}(\cV)$ is the category of coloured symmetric operads in $\cV$.
\end{example}

\begin{example}
If we consider the $\FF$-multicategory $\und{\NN}$ of \ref{exnonsym}, then the category $\operc{\und{\NN}}(\cV)$ is the category of multicategories enriched in $\cV$.
\end{example}

\begin{construction}\label{funcop}
For every multifunctor $\phi \colon \MM \to \NN$, we define a natural transformation $(-)^*\phi \colon (-)^*\MM \circ (\phi_0)^* \Rightarrow (-)^*\NN$ with $(\phi_0)^*\colon \Psh(\NN_0) \to \Psh(\MM_0)$ the functor obtained by pulling back along $\phi_0$.
For every functor $F\colon \cA \to \NN_0$, if we denote by $G$ the functor $(\phi_0)^*F$, we have to construct a multifunctor
$F^*\phi \colon G^* \MM \to F^* \NN$. At the level of categories of objects, the functor $(F^*\phi )_0 \colon \MM_0 \times_{\NN_0}\cA \to  \cA$ is the canonical projection and at the level of multimorphisms, the functor 
$$
(F^*\phi )_1 \colon \MM_1 \times_{\left( \MM_0 \times \FF \MM_0 \right)}\left(\left(\cA\times_{\NN_0}\MM_0 \right) \times \FF \left(\cA \times_{\NN_0}\MM_0 \right) \right) \to \NN_1 \times_{\left( \NN_0 \times \FF \NN_0 \right)}\left( \cA \times \FF \cA \right)
$$
in given by $\phi_1 \times_{\Id{(\FF \NN_0 \times \NN_0)}} \Id{\FF \cA \times \cA}$ using the equivalence
$$
\MM_1 \times_{\left( \MM_0 \times \FF \MM_0 \right)}\left(\left(\cA\times_{\NN_0}\MM_0 \right) \times \FF \left(\cA \times_{\NN_0}\MM_0 \right) \right)\simeq \MM_1 \times_{\left( \NN_0 \times \FF \NN_0 \right)}(\FF \cA \times \cA).
$$
In particular, for every multifunctor $\phi \colon \MM \to \NN$, we obtain a natural transformation $\Phi_\phi\colon \Phi_\NN \Rightarrow \Phi_{\MM}\circ (\phi_0)^*$ of functors $\opo{\Psh(\NN_0)} \to \CAT$ by precomposing and post-composing $(-)^*\phi$ by $\el$ and $\oper{(-)}(\cV)$ respectively. We thus obtain a fuctor $\phi^* \colon \operc{\NN}(\cV) \to \operc{\MM}(\cV)$ on the associated Grothendieck constructions using \cite[Section 2.3]{grothmod}.
\end{construction}

\begin{proposition}\label{exundcat}
The functor $(-)_0 \colon \multi{\FF} \to \Cat$ that sends an $\FF$-multicategory $\MM$ to its category of objects $\MM_0$ admits a left adjoint that sends a category $\cC$ to the $\FF$-multicategory $\widetilde{\cC}$ constant on $\cC$.
For every presheaf $A \colon \opo{\cC} \to \Set$, the category of $A$-coloured $\widetilde{\cC}$-operads in $\cV$ is the category $\fun(\opo{\cC},\Cat_{\cV})^A$ of functors $F\colon \opo{\cC} \to \Cat_{\cV}$ (with $\Cat_{\cV}$ the category of small $\cV$-categories) such that the composition $\opo{\cC} \overset{F}{\longrightarrow} \Cat_{\cV} \overset{\operatorname{Ob}}{\longrightarrow} \Set$ is $A$. Moreover, the category of coloured $\widetilde{\cC}$-operads is the category $\fun (\opo{\cC},\Cat_{\cV})$.
\begin{proof}
Let $A \colon \opo{\cC} \to \Set$ be a presheaf and let $\cO$ be a $A$-coloured $\widetilde{\cC}$-operad in $\cV$. The category $\left(\widetilde{\cC}^A\right)_1$ is the pullback $\el (X) \times_{\cC} \el (A)$. For every object $X$ of $\cC$, we can define a $\cV$-category $\cO(X)$ with set of objects $A(X)$ and, given two objects $x,y$ in $A(X)$, whose object in $\cV$ of morphisms is $\cO(X,x,y)$ where $(X,x,y)$ is seen as an object of $\el (A) \times_{\cC} \el (A)$. The composition and identities of the operad $\cO$ respectively provide a composition and identities for the $\cV$-category $\cO(X)$. This construction provides the promised equivalence between $A$-coloured $\widetilde{\cC}$-operads in $\cV$ and $\fun(\opo{\cC},\Cat_{\cV})^A$. The last part of the proposition follows from the fact that the category $\fun(\opo{\cC},\Cat_{\cV})$ is the Grothendieck construction of the functor $\opo{\cC} \to \CAT$ that sends $A$ to $\fun(\opo{\cC},\Cat_{\cV})^A$.
\end{proof}
\end{proposition}

\begin{definition}\label{inducedcat}
For every $\FF$-multicategory $\MM$, if we denote by $i\colon \widetilde{\MM_0} \to \MM$
the counit of the adjunction $\widetilde{(-)}\colon \Cat \rightleftarrows \multi{\FF} \colon (-)_0$ of \ref{exundcat}, using \ref{funcop}, we obtain a functor $i^* \colon \operc{\MM}(\cV) \to \fun (\opo{\MM_0},\Cat_{\cV})$. For every coloured $\MM$-operad $\cO$ in $\cV$, the functor $i^*(\cO)\colon \opo{\MM_0} \to \Cat_{\cV}$ is the \emph{underlying parametrized $\cV$-category of $\cO$}.
\end{definition}

\subsection{Parametrized operads}

In this section, we use the formalism of $\MM$-operads to define parametrized operads. We recall first the notion of orbital subcategories that generalizes Nardin's notion of orbital categories \cite[Definition 4.1]{dede}.

\begin{definition}[\protect{\cite[Definition 4.2.2]{CLL1}}]\label{deforb}
Let $\sO$ be a small category. An \emph{orbital subcategory $\sT$ of $\sO$} is a wide subcategory $\sT$ of $\sO$ such that the pullback in $\FF \sO$ of any morphism in $\FF\sT\subseteq \FF \sO$ exists and belongs to $\FF \sT$. An \emph{orbital pair} is a pair $(\sO,\sT)$ such that $\sO$ is a small category and such that $\sT$ is an orbital subcategory of $\sO$.
A \emph{functor $F\colon (\sO_1,\sT_1) \to (\sO_2,\sT_2)$ of orbital pairs} is a functor $F\colon \sO_1 \to \sO_2$ such that $F(\sT_1)\subseteq \sT_2$ and such that $\FF F$ preserves pullbacks of the form
$$
\begin{tikzcd}
A \arrow[r] \arrow[d] & B \arrow[d,"g"] \\
C \arrow[r] & D
\end{tikzcd}
$$
where $g$ belongs to $\sT_1$. A category $\sO$ is \emph{orbital} if $\sO$ is an orbital subcategory of itself.
\end{definition}

\begin{proposition}[\protect{\cite[Observation A.5]{nat1}}]\label{transorb}
For every finite group $G$, the poset of orbital subcategories $\suborb{\sO_G}$ of the orbit category $\sO_{G}$ is equivalent to the poset $\operatorname{Tr}(G)$ of transfer systems for $G$ is the sense of \cite[Lemma 6]{BBR}. In particular, it follows from \cite{BH1}, \cite{Rub1}, \cite{BonPer1} and \cite{BBR}  that the poset $\suborb{\sO_G}$ is also isomorphic to the poset of $N_{\infty}$-operads and of $G$-indexing systems. 
\end{proposition}

\begin{example}
For every category $\cC$, the core of $\cC$ is an orbital subcategory of $\cC$. We thus obtain an orbital pair $(\cC,\cC^{\simeq})$. 
\end{example}

\begin{example}
The category of finite sets of cardinality less than $n$ and surjective functions is an orbital category.
\end{example}

\begin{definition}\label{defparaop}
For every orbital pair $(\sO,\sT)$, consider the subcategory $\Sigma_{(\sO,\sT)}$ of $\Arr(\FF \sO)$ whose objects are morphisms $f\colon A \to B$ in $\FF\sT$ such that $B$ belongs to $\sO$ and where a morphism from $f\colon A \to B$ to $g\colon C \to D$ is a pullback square
$$
\begin{tikzcd}
A \arrow[r,"i"] \arrow[d,"f"'] & C \arrow[d,"g"] \\
B \arrow[r,"h"'] & D.
\end{tikzcd}
$$
in $\FF \sO$. We now define an $\FF$-multicategory $\und{\Sigma_{(\sO,\sT)}}$ whose category of multimorphisms is $\Sigma_{(\sO,\sT)}$, whose category of objects is $\sO$, whose identity functor $1_{\und{\Sigma_{(\sO,\sT)}}}\colon \sO \to \Sigma_{(\sO,\sT)}$ sends an object $X$ to the identity $\Id{X}$, and whose composition functor $\odot \colon \Sigma_{(\sO,\sT)}\times_{\FF \sS} \FF \Sigma_{(\sO,\sT)}$ sends $(f,(f_1,\ldots,f_n))$ to the composition $f\circ (f_1\sqcup \cdots \sqcup f_n)$ in the category $\FF \sT$.
A (coloured) \emph{$(\sO,\sT)$-operad in $\cV$} is a (coloured) $\und{\Sigma_{(\sO,\sT)}}$-operad in $\cV$ and we denote by $\oper{(\sO,\sT)}(\cV)$ (respectively by $\operc{(\sO,\sT)}(\cV)$ ) the category of $(\sO,\sT)$-operads (respectively coloured $(\sO,\sT)$-operads) in $\cV$. If $\sO$ is an orbital category, then we denote by $\oper{\sO}(\cV)$ (respectively by $\operc{\sO}(\cV)$) the category of $(\sO,\sO)$-operads in $\cV$ (respectively coloured $(\sO,\sO)$-operads in $\cV$).
\end{definition}

\begin{remark}
If $\sO$ is an orbital category, a coloured $(\sO,\sO)$-operad in $\sSet$ is a simplicial coloured $\sO$-operad in the sense of \cite[Definition 2.5.4]{DJ}. 
\end{remark}

\begin{remark}
Note that it should be possible to extend \ref{defparaop} to the setting of weak indexing systems in the sense of \cite{nat1}.
\end{remark}

\begin{proposition}\label{funcorb}
Every functor of orbital pairs $F\colon (\sO_1,\sT_1) \to (\sO_2,\sT_2)$ induces a multifunctor $F\colon \und{\Sigma_{(\sO_1,\sT_1)}} \to \und{\Sigma_{(\sO_2,\sT_2)}}$. In particular, every functor of orbital pairs $F\colon (\sO_1,\sT_1) \to (\sO_2,\sT_2)$ induces by restriction a functor $F^* \colon \oper{(\sO_1,\sT_1)}(\cV) \to \oper{(\sO_2,\sT_2)}(\cV)$ as well as a functor 
$F^* \colon \operc{(\sO_1,\sT_1)}(\cV) \to \operc{(\sO_2,\sT_2)}(\cV)$.
\begin{proof}
The functor $F$ induces a functor $\Arr(\FF F)\colon \Arr(\FF \sO_1) \to \Arr(\FF \sO_2)$ and because $F \cT_1 \subseteq \cT_2$ and because $\FF F$ preserves pullbacks of morphisms in $\cT$, we deduce that  $\Arr(\FF F)$ restricts to a functor $\Sigma_{F}\colon \Sigma_{(\sO_1,\sT_2)} \to \Sigma_{(\sO_2,\sT_2)}$. Thus, we obtain a multifunctor $\und{\Sigma_F}\colon \und{\Sigma_{(\sO_1,\sT_2)}} \to \und{\Sigma_{(\sO_2,\sT_2)}}$. The proposition now follows from \ref{funcop}.
\end{proof}
\end{proposition}

\begin{example}
The final category $*$ is an orbital category and the $\FF$-multicategory $\und{\Sigma_*}$ is the $\FF$-multicategory that we already defined in \ref{exsymmop}. In particular, the category $\operc{*}(\cV)$, is the category $\operc{}(\cV)$ of coloured symmetric operads in $\cV$.
\end{example}

\begin{example}
For every category $\cC$, the category $\operc{(\cC,\cC^{\simeq})}(\cV)$ is the category $\fun(\opo{\cC},\operc{}(\cV))$. In particular, for every discrete group $G$, $\operatorname{B}G$ is an orbital category and the category of $\operatorname{B}G$-operads in $\cV$ is the category $\operc{}(\cV) ^G$ of $G$-objects in $\operc{}(\cV)$.
\end{example}

\begin{example}\label{exgenop}
For every $G$-transfer system $\cT$, if follows from \ref{transorb} that we obtain a notion of \emph{$\cT$-incomplete coloured operads in $\cV$} given by the category $\operc{(\sO_G,\sO_{\cT})}(\cV)$. This category is expected to be equivalent to the category of \emph{partial genuine equivariant operads} of Bonventre and Pereira \cite{BonPer1}. However, such a comparison has not yet been made.
\end{example}

\begin{proposition}\label{elmop}
For every $G$-transfer system $\cT$, the inclusion $\operatorname{B}G\hookrightarrow \sO_G$ produces a functor $(\operatorname{B}G,\operatorname{B}G) \to (\sO_G,\sO_{\cT})$ of orbital pairs. In particular, using \ref{funcorb}, we obtain a restriction functor $\operc{(\sO_G,\sO_{\cT})} \to\operc{}(\cV)^G$.
\begin{proof}
Note first that the groupoid $\operatorname{B}G$ is contained in every $\sO_{\cT}$ for every $G$-transfer system $\cT$. Moreover, it is clear that the functor $\FF \operatorname{B}G \hookrightarrow \FF \sO_G$ preserves pullbacks. It follows that for every $G$-transfer system $\cT$, the inclusion $\operatorname{B}G \hookrightarrow \sO_G$ provides a functor of orbital pairs $(\operatorname{B}G,\operatorname{B}G) \to (\sO_G,\sO_{\cT})$.
\end{proof}
\end{proposition}

\section{Convolution product for double categories}\label{sec3}
The goal of this section is to construct a convolution product for double categories that generalizes Day's convolution product. Note that similar considerations have already been made recently in \cite{convdoub}.
Until the end of this section, $\DD$ is a small double category, $(\cC,\otimes_{\cC})$ is a cocomplete closed monoidal category, and $\cV$ is a cocomplete cartesian closed category.

\subsection{Monoidal structure on spans}

\begin{definition}[\protect{\cite[Definition 4.3]{BonPer1}}]
If $\cB$ and $\cC$ are categories, then we denote by $\wspan(\cB,\cC)$ the category with objects the spans of functors
$$
\cB \overset{k}{\longleftarrow} \cA \overset{X}{\longrightarrow} \cC
$$
where $\cA$ is a small category, with morphisms the following diagrams
$$
\begin{tikzcd}[row sep=small]
& \cA \arrow[rd,"X"{name=U}]  \arrow[ld,"k"'] \arrow[dd,"i"] & \\
\cB & & \cC \\
& \cD \arrow[ur,"Y"'] \arrow[Rightarrow,"\varphi",shorten >=3mm,shorten <=3mm,from=U] \arrow[ul,"l"] &
\end{tikzcd}
$$
which we write as $(i,\varphi)\colon (k_1,X_1)\to (k_2,X_2)$, and composition given in the obvious way.
\end{definition}

\begin{proposition}
The category $\wspan (\DD_1, \cC)$ is endowed with the structure of a monoidal category defined as follows: the tensor product $\otimes$ of two spans $\DD_1 \overset{k_1}{\longleftarrow} \cA_1 \overset{X_1}{\longrightarrow} \cC$ and $\DD_1 \overset{k_2}{\longleftarrow} \cA_2 \overset{X_2}{\longrightarrow} \cC$ is the span
$$
\DD_1 \overset{k}{\longleftarrow} \cA_1 \times_{\DD_0} \cA_2 \overset{X}{\longrightarrow} \cC
$$
where the left hand morphism is the composition $\cA_1 \times_{\DD_0} \cA_2 \overset{k_1\times_{\DD_0} k_2}{\longrightarrow} \DD_1 \times_{\DD_0} \DD_1 \overset{\odot}{\longrightarrow} \DD_1$ and where the right hand side morphism is the composition $\cA_1 \times_{\DD_0} \cA_2 \overset{X_1\times X_2}{\longrightarrow} \cC \times \cC \overset{\otimes_{\cC}}{\longrightarrow} \cC$. The tensor product $\otimes$ of a pair of morphisms of spans 
$$
\begin{tikzcd}[row sep=small]
& \cA_1 \arrow[rd,"X_1"{name=U}]  \arrow[ld,"k_1"'] \arrow[dd,"i_1"'] & & & \cA_2 \arrow[rd,"X_2"{name=V}]  \arrow[ld,"k_2"'] \arrow[dd,"i_2"'] &\\
\DD_1 & & \cC & \DD_1 & & \cC \\
& \cB_1 \arrow[ur,"Y_1"'] \arrow[Rightarrow,"\varphi_1"',shorten >=3mm,shorten <=3mm,from=U] \arrow[ul,"l_1"] & & & \cB_2 \arrow[ur,"Y_2"'] \arrow[Rightarrow,"\varphi_2"',shorten >=3mm,shorten <=3mm,from=V] \arrow[ul,"l_2"] &
\end{tikzcd}
$$
is the morphisms of spans $(i_1 \times i_2, \varphi_1 \times \varphi_2)$ obtained by considering the diagram
$$
\begin{tikzcd}[row sep=small,column sep=large]
& & \cA_1\times_{\DD_0} \cA_2 \arrow[dd,"i_1\times i_2"'] \arrow[r] \arrow[ld,"k_1\times_{\DD_0}k_2"'] & \cA_1 \times \cA_2  \arrow[dd,"i_1\times i_2"'] \arrow[rd,"X_1\times X_2"{name=U}] & &   \\
\DD_1& \arrow[l,"\odot"'] \DD_1 \times_{\DD_0} \DD_1 & & & \cC \times \cC \arrow[r,"\otimes_{\cC}"] & \cC.\\
& & \cB_1\times_{\DD_0}\cB_2   \arrow[ul,"l_1\times_{\DD_0}l_2"] \arrow[r] &  \cB_1 \times \cB_2   \arrow[ur,"Y_1\times Y_2"'] 
 \arrow[Rightarrow,"\varphi_1 \times \varphi_2"',shorten >=3mm,shorten <=3mm,from=U] & &
\end{tikzcd}
$$
The unit of $\wspan(\DD_1,\cC)$ is the span 
$$
\DD_1 \overset{1_{\MM}}{\longleftarrow} \DD_0 \overset{1_{\cC}}{\longrightarrow}\cC
$$
where the right hand side functor is constant on the unit $1_{\cC}$ of $\cC$.
\begin{proof}
We have to construct the associators and unitors of the tensor product $\otimes$.
If we consider three spans $S_1=\DD_1 \overset{k_1}{\longleftarrow} \cA_1 \overset{X_1}{\longrightarrow} \cC$, $S_2=\DD_1 \overset{k_2}{\longleftarrow} \cA_2 \overset{X_2}{\longrightarrow} \cC$ and $S_3=\DD_1 \overset{k_3}{\longleftarrow} \cA_3 \overset{X_3}{\longrightarrow} \cC$, their tensor product $S_1 \otimes (S_2 \otimes S_3)$ is the span obtained as the composition of
$$
\DD_1 \overset{\odot}{\longleftarrow} \DD_1\times(\DD_1 \times_{\DD_0}\DD_1) \overset{}{\longleftarrow} \cA_1 \times_{\DD_0}(\cA_2 \times_{\DD_0}\cA_3)\overset{}{\longrightarrow} \cC \times (\cC \times \cC) \overset{\otimes}{\longrightarrow} \cC
$$
and the iterated tensor product $(S_1\otimes S_2)\otimes S_3$ is also obtained in a similar way. It follows that the associator $\alpha_{S_1,S_2,S_3}\colon S_1\otimes (S_2 \otimes S_3) \simeq (S_1 \otimes S_2)\otimes S_3$ can be provided by the following diagram 
$$
\begin{tikzcd}[column sep=large]
&\DD_1 \times_{\DD_0}(\DD_1 \times_{\DD_0}) \DD_1 \arrow[dd,"\simeq"] \arrow[ld,"\odot"'] & \cA_1\times_{\DD_0} (\cA_2\times_{\DD_0}\cA_3) \arrow[dd,"\simeq"'] \arrow[r,"X_1 \times (X_2 \times X_3)"] \arrow[l,"k_1\times_{\DD_0} (k_2\times_{\DD_0}k_3)"'] & \cC \times (\cC \times \cC) \arrow[dd,"\simeq"'] \arrow[rd,"\otimes_{\cC}"{name=U}] &    \\
\DD_1&  & & & \cC   \\
&(\DD_1 \times_{\DD_0}\DD_1) \times_{\DD_0} \DD_1 \arrow[lu,"\odot"] & (\cA_1\times_{\DD_0}\cA_2)\times_{\DD_0}\cA_3   \arrow[l,"(k_1\times_{\DD_0}k_2)\times_{\DD_0}k_3"] \arrow[r,"(X_1 \times X_2)\times X_3"'] &  (\cC \times \cC)\times \cC   \arrow[ur,"\otimes_{\cC}"'] 
 \arrow[Rightarrow,"\beta"',shorten >=3mm,shorten <=3mm,from=U] & 
\end{tikzcd}
$$
where $\beta$ is the associator of $\cC$. The pentagon axiom can be verified using the observation that it is verified both in $\cC$ and for pullbacks over a given category. For every span $S=\DD_1 \overset{k}{\longleftarrow} \cA \overset{X}{\longrightarrow} \cC$, the right unitor $\eta_{S}\colon  S\otimes 1\simeq S$ is given by the diagram
$$
\begin{tikzcd}[row sep=small]
& \cA\times_{\DD_0} \DD_0  \arrow[ld,""'] \arrow[dd,"\simeq"'] \arrow[r] & \cC\times * \arrow[dd,"\simeq"'] \arrow[r] & \cC \times \cC \arrow[rd,"\otimes_{\cC}"] \arrow[ldd,Rightarrow,shorten >=4mm,shorten <=4mm,"\rho"] \arrow[dd,"\otimes_{\cC}"] &  \\
\DD_1 &   & & & \cC \\
& \cA \arrow[r,"X"'] \arrow[ul,"k"] &\cC \arrow[r,"\Id{\cC}"'] & \cC \arrow[ru,"\Id{\cC}"'] &
\end{tikzcd}
$$
where $\rho$ is the right unitor of $\cC$. The left unitor is given by a similar diagram, and the triangle identities are easily verified using that they are already verified in $\cC$ and for pullbacks over a given category.
\end{proof}
\end{proposition}

\begin{remark}
Note that the tensor product $\otimes$ on $\wspan (\DD_1, \cC)$ exists only if $\DD$ is a strict double category.
\end{remark}

\subsection{Left Kan extension along morphisms of cocartesian fibrations}
The results of this subsection are mostly technical and will be used in \ref{secconv}.
\begin{notation}
Consider a functor $G\colon \cA \to \cB$ over $\cC$
$$
\begin{tikzcd}[column sep=small]
\cA \arrow[dr,"F"'] \arrow[rr,"G"] && \cB \arrow[dl,"H"] \\
&\cC.&
\end{tikzcd}
$$
Recall that, for $X$ in $\cB$, the category $G_{/X}$ is obtained as the following pullback
$$
\begin{tikzcd}
G_{/X} \arrow[d] \arrow[r] & \cB_{/X} \arrow[d] \\
\cA \arrow[r] & \cB.
\end{tikzcd}
$$
For every object $X$ of $\cB$, we can consider the fiber category $H^{-1}(HX) \hookrightarrow \cB$ and we thus obtain a subcategory of $\left(H^{-1}(HX)\right)_{/X} \hookrightarrow \cB_{/X}$. We denote by $G_{/_{H}X}$ the subcategory of $G_{/X}$ obtained as the pullback 
$$
\begin{tikzcd}
G_{/_{H}X} \arrow[d,hook] \arrow[r] & \left(H^{-1}(HX)\right)_{/X} \arrow[d,hook] \\
G_{/X} \arrow[d] \arrow[r] & \cB_{/X} \arrow[d] \\
\cA \arrow[r,"G"'] & \cB.
\end{tikzcd}
$$
    
\end{notation}

\begin{proposition}[\protect{\cite[Lemma 3.1.11]{Haug1}}]\label{refsubt}
If $G\colon \cA \to \cB$ is a morphism of cocartesian fibrations
$$
\begin{tikzcd}[column sep=small]
\cA \arrow[dr,"F"'] \arrow[rr,"G"] && \cB \arrow[dl,"H"] \\
&\cC,&
\end{tikzcd}
$$
then the subcategory $G_{/_H X}$ is a reflexive subcategory of $G_{/X}$. In particular, the inclusion $G_{/_H X} \hookrightarrow G_{/X}$ is final.
\end{proposition}

\begin{corollary}\label{fiberfi}
Consider a morphism of left fibrations 
$$
\begin{tikzcd}[column sep=small]
\cA \arrow[dr,"F"'] \arrow[rr,"G"] && \cB \arrow[dl,"H"] \\
&\cC&
\end{tikzcd}
$$
such that $G$ is an isofibration. For every object $X$ of $\cB$, the inclusion $G^{-1}(X) \hookrightarrow G_{/X}$ is a final functor.
\begin{proof}
The functor $H$ is a left fibration, and it follows that $H^{-1}(H X)$ is a groupoid. In particular, the category $\left(H^{-1}(HX)\right)_{/X}$ is contractible. We now obtain the diagram
$$
\begin{tikzcd}
G^{-1}(X) \arrow[r] \arrow[d,hook] & * \arrow[d,"\simeq"]\\
G_{/_{H}X} \arrow[d] \arrow[r] & \left(H^{-1}(HX)\right)_{/X} \arrow[d] \\
G_{/X}  \arrow[r] & \cB_{/X}.  \\
\end{tikzcd}
$$
Finally, because $G$ is an isofibration, the functor $G_{/_{H}X} \to \left(H^{-1}(HX)\right)_{/X}$ is also an isofibration and $G^{-1}(X) \to G_{/_{H}X}$ is an equivalence as it is the pullback of an equivalence along an isofibration.
\end{proof}
\end{corollary}

\subsection{Convolution product for double categories}\label{secconv}

\begin{definition}
A double category $\DD$ is \emph{target left fibrant} if the target functor $t\colon \DD_1 \to \DD_0$ is a left fibration and if the horizontal composition functor $\odot \colon \DD_1\times_{\DD_0}\DD_1 \to \DD_1$ is an isofibration. If $\DD$ is a target left fibrant double category, then a functor $F\colon \cA \to \DD_1$ is a \emph{target left fibration} if it is an isofibration and if the composition $\cA \overset{F}{\longrightarrow}\DD_1 \overset{t}{\longrightarrow} \DD_0$ is a left fibration.
\end{definition}

For the rest of this section, we suppose that $\DD$ is a small target left fibrant double category.

\begin{definition}
The category $\tspan (\DD_1,\cC)$ is the full subcategory of $\wspan (\DD_1,\cC)$ spanned by the spans $\DD_1 \overset{k}{\longleftarrow} \cA \overset{X}{\longrightarrow} \cC$ such that $k$ is a target left fibration.  
\end{definition}

\begin{proposition}
The tensor product $\otimes$ on $\wspan (\DD_1,\cC)$ restricts to a monoidal structure on the subcategory $\tspan (\DD_1,\cC)$.
\begin{proof}
Consider two spans  $\DD_1 \overset{k}{\longleftarrow} \cA \overset{X}{\longrightarrow} \cC$ and  $\DD_1 \overset{l}{\longleftarrow} \cB \overset{Y}{\longrightarrow} \cC$ such that $k$ and $l$ are target left fibrations. The left leg of their tensor product is the following composition 
$$
\cA \times_{\DD_0} \cB \overset{k\times_{\DD_0}l}{\longrightarrow} \DD_1 \times_{\DD_0} \DD_1 \overset{\odot}{\longrightarrow} \DD_1
$$
and its composition with $t\colon \DD_1 \to \DD_0$ is the composition
$$
\cA \times_{\DD_0} \cB \overset{p_2}{\longrightarrow} \cB \overset{l} {\longrightarrow} \DD_1 \overset{t}{\longrightarrow} \DD_0
$$
that is a left fibration as a composition of left fibrations (using that $\cA \times_{\DD_0} \cB \to \cB $ is a pullback of $t\circ k$).
\end{proof}
\end{proposition}

\begin{proposition}\label{LKA}
The functor 
$$
\Lan \colon \wspan (\DD_1,\cC) \to \fun (\DD_1, \cC)
$$
that sends a span $\DD_1 \overset{k}{\longleftarrow} \cA \overset{X}{\longrightarrow} \cC$ to the left Kan extension of $X$ along $k$ is a left adjoint of the functor $v\colon \fun (\DD_1,\cC) \to \wspan (\DD_1,\cC)$ that sends a functor $X\colon \DD_1 \to \cC$ to the span $\DD_1 \overset{\Id{\DD_1}}{\longleftarrow} \DD_1 \overset{X}{\longrightarrow} \cC$.
Moreover, this adjunction restricts to the adjunction 
$$
\Lan \colon \tspan(\DD_1,\cC) \rightleftarrows \fun (\DD_1,\cC)\colon v.
$$
\begin{proof}
The first part of the proposition is \cite[Remark 4.5]{BonPer1}. The second part of the statement follows from the fact that the identity of $\DD_1$ is a target left fibration.
\end{proof}
\end{proposition}

\begin{proposition}[\protect{\cite[Theorem 3.1]{Kel2}}]\label{indmon}
For every reflexive adjunction $L\colon \cD \rightleftarrows \cE \colon R$ such that the category $\cD$ is endowed with the structure of a monoidal category $(\cD,\otimes_{\cD})$, if for all $A$ and $B$ in $\cD$, the morphism
$$
L(\eta_A \otimes_{\cD} \eta_B)\colon L(A\otimes B) \to L\left(RL A \otimes_{\cD} RL B \right)
$$
is an isomorphism, then $\cE$ is naturally endowed with the structure of a monoidal category where the tensor product of two objects $X$ and $Y$ of $\cE$ is defined by $X\otimes_{\cE} Y :=L(R(X)\otimes_{\cD}R(Y))$ and where the unit of $\cE$ is the image of the unit of $\cD$ under $L$. Moreover, the category of algebras in $\cE$ over $\otimes_{\cE}$ is the full subcategory of the category of algebras in $\cD$ over $\otimes_{\cD}$ spanned by the algebras whose underlying object belongs to $\cE$. Furthermore, if $X_1,\ldots ,X_n$ are objects of $\cE$, then we have $X_1 \otimes_{\cE} \cdots \otimes_{\cE} X_n\simeq L(R(X)\otimes_{\cD} \cdots \otimes_{\cD} R(X_n))$.
\end{proposition}

\begin{proposition}\label{tensleftkan}
 For all functors $F\colon \cA_1 \to \cC$, $ G \colon \cB_1 \to \cC$, $H\colon \cA_1 \to \cB_1$ and $L\colon \cA_2 \to \cB_2$ such that $\cA_1$, $\cA_2$, $\cB_1$ and $\cB_2$ are small categories, if we denote by $\operatorname{L}_H F$ and $\operatorname{L}_L G$ the left Kan extension of $F$ along $H$ and $G$ along $L$ respectively, then the left Kan extension of $\otimes_{\cC} \circ (F\times G)$ along $H \times L$ is $\otimes_{\cC} \circ (\operatorname{L}_{H}F \times \operatorname{L}_L G)$.
\begin{proof}
This follows directly from the pointwise formula for left Kan extensions using that $\otimes_{\cC}$ preserves colimits in both variables.
\end{proof}
\end{proposition}

\begin{proposition}\label{convind}
The adjunction
$$
\Lan \colon \tspan (\DD_1,\cC) \rightleftarrows \fun(\DD_1, \cC)\colon v
$$
satisfies the conditions of \ref{indmon}.
\begin{proof}
Consider two spans $\DD_1 \overset{F}{\longleftarrow} \cA \overset{X}{\longrightarrow} \cC$ and  $\DD_1 \overset{G}{\longleftarrow} \cB \overset{Y}{\longrightarrow} \cC$ such that $F$ and $G$ are target left fibrations, and consider their associated left Kan extensions 
$$
\begin{tikzcd}
\cA \arrow[r,"X"{name=U}]   \arrow[d,"F"'] & \cC  & \cB \arrow[r,"Y"{name=V}]   \arrow[d,"G"'] & \cC \\
 \DD_1 \arrow[ur,"\operatorname{L}_F X"'] \arrow[Rightarrow,"\varphi"',shorten >=3mm,shorten <=3mm,from=U]  &  & \DD_1. \arrow[ur,"\operatorname{L}_G Y"'] \arrow[Rightarrow,"\psi"',shorten >=3mm,shorten <=3mm,from=V]  &
\end{tikzcd}
$$
We have to show that the composite of the following diagram 

\begin{equation}\tag{E}\label{convdiag}
\begin{tikzcd}[row sep=small]
 \cA\times_{\DD_0} \cB \arrow[dd,"F\times_{\DD_0} G"'] \arrow[r]  & \cA \times \cB  \arrow[dd,"F\times G"'] \arrow[rd,"X \times Y"{name=U}] & &   \\
 & & \cC \times \cC \arrow[r,"\otimes_{\cC}"] & \cC\\
 \DD_1 \times_{\DD_0}\DD_1    \arrow[r] &  \DD_1 \times \DD_1   \arrow[ur,"\left(\operatorname{L}_F X\right)\times \left(\operatorname{L}_G Y\right)"'] 
 \arrow[Rightarrow,"\varphi\times \psi"',shorten >=3mm,shorten <=3mm,from=U] & &
\end{tikzcd}
\end{equation}

is a left Kan extension diagram. 
Using \ref{tensleftkan} and the point-wise formula for left Kan extension, it is enough to show that the functor 
$$
\Phi \colon \left(F\times_{\DD_0} G\right)_{/(X,Y)} \to \left(F\times G \right)_{/(X,Y)}
$$
is final.
The left hand side square of the diagram \ref{convdiag} is a pullback, and we thus obtain the following square

\begin{equation}\label{sqinclu}\tag{F}
\begin{tikzcd}
\left(F\times_{\DD_0} G\right)^{-1}(X,Y) \arrow[d,hook] \arrow[r,"\simeq "]  & (F\times G)^{-1}(X,Y) \arrow[d,hook] \\
\left(F\times_{\DD_0} G\right)_{/(X,Y)} \arrow[r,"\Phi"']   &   \left(F\times G \right)_{/(X,Y)}. \\
\end{tikzcd}
\end{equation}

The functor $F\times G$ is an isofibration as a product of isofibrations, and it follows that $F\times_{\DD_0}G$ is also an isofibration as the pullback of an isofibration. Moreover, $F\times G$ can be seen as a morphism of left fibrations as follows
$$
\begin{tikzcd}[column sep=tiny]
\cA \times \cB \arrow[rd,"(t\circ F)\times (t\circ G)"'] \arrow[rr,"F\times G"] & & \DD_1 \times \DD_1 \arrow[ld,"t\times t"] \\
&\DD_0 \times \DD_0 &
\end{tikzcd}
$$
and the functor $F\times_{\DD_0}G$ can be seen as a morphism of left fibrations as follows
$$
\begin{tikzcd}[column sep=tiny]
\cA\times_{\DD_0} \cB \arrow[rr,"F\times_{\DD_0} G"] \arrow[rd,"t\circ G \circ p_2"'] & & \DD_1 \times_{\DD_0} \DD_1 \arrow[dl,"t\circ p_2"] \\
&\DD_0. &
\end{tikzcd}
$$
Finally, by \ref{fiberfi}, the vertical inclusions of the square \ref{sqinclu} are final functors, and it follows that $\Phi$ is also final.
\end{proof}
\end{proposition}

\begin{corollary}\label{algmon}
The category $\fun (\DD_1,\cC)$ is naturally endowed with the structure of a closed monoidal category for which the tensor product $\otimes_{\DD}$ of two functors $F,G \colon \DD_1 \to \cC$ is given by the following left Kan extension 
$$
\begin{tikzcd}
\DD_1 \times_{\DD_0} \DD_1 \arrow[r,"F\times G"] \arrow[d,"\odot "'] & \arrow[dl,Rightarrow,shorten >=4mm,shorten <=4mm] \cC \times \cC \arrow[r,"\otimes_{\cC}"] & \cC, \\
\DD_1 \arrow[rru,bend right=10,"F\otimes_{\DD} G"'] &&
\end{tikzcd}
$$
and whose unit is the left Kan extension
$$
\begin{tikzcd}
\DD_0 \arrow[d,"1_{\DD}"'] \arrow[r] &* \arrow[ld,Rightarrow,shorten >=3mm,shorten <=3mm] \arrow[r,"1_{\cC}"]&  \cC. \\
\DD_1 \arrow[rru,bend right=10] &&
\end{tikzcd}
$$
Moreover, the category of algebras over $\otimes_{\DD}$ is the category $\lax(\DD,\cC)$. Furthermore, for all functors $F_1,\ldots, F_n  \colon \DD_1 \to \cC$, the tensor product $F_1 \otimes_{\DD} \cdots \otimes_{\DD} F_n$ is the left Kan extension
$$
\begin{tikzcd}
\DD_n \arrow[d,"\odot^n"'] \arrow[r] & \cC^{\times n} \arrow[ld,Rightarrow,shorten >=3mm,shorten <=3mm] \arrow[r,"\otimes_{\cC}^n"] & \cC \\
\DD_1 \arrow[rru,bend right=10,"F_1 \otimes_{\DD} \cdots \otimes_{\DD} F_n"'] && 
\end{tikzcd}
$$
where $\DD_n$ is the pullback $\DD_1 \times_{\DD_0} \cdots \times_{\DD_0} \DD_1$ with $n$ occurences of $\DD_1$.
\begin{proof}
The existence of the monoidal structure follows from the combination of \ref{convind} and \ref{indmon}. We show now that $(\fun (\DD_1,\cC),\otimes_{\DD})$ is a closed monoidal category. We will show that $(\fun (\DD_1,\cC),\otimes_{\DD})$ is left closed, the right closedness can be shown with a symmetric proof.
Let $F$ be a functor $F\colon \DD_1 \to \cC$. We have to show that the functor 
$$
(-)\otimes_{\DD} F \colon \fun(\DD_1,\cC) \to \fun(\DD_1,\cC)
$$
admits a right adjoint. This functor can be viewed as the following composition
$$
\fun\left(\DD_1,\cC\right) \overset{\otimes_{\cC}\circ((-)\times F)}{\longrightarrow} \fun\left(\DD_1 \times \DD_1,\cC\right) \overset{}{\longrightarrow} \fun\left(\DD_1 \times_{\DD_0} \DD_1,\cC\right) \overset{\odot_!}{\longrightarrow} \fun(\DD_1,\cC).
$$
The last functor is a left adjoint of $\odot^*$ and the second functor admits a right adjoint given by taking right Kan extensions along the functor $\DD_1 \times_{\DD_0} \DD_1 \to \DD_1 \times \DD_1$. It remains to show that the first functor admits a right adjoint. As the monoidal category $\cC$ is closed, this functor is a left adjoint of the functor $\fun (\DD_1 \times \DD_1,\cC) \to \fun (\DD_1,\cC)$ that sends a functor $H\colon \DD_1 \times \DD_1  \to \cC$ to the functor $\DD_1 \to \cC$ that sends an object $X$ of $\DD_1$ to the end
$$
\int_{Y\in \DD_1} \hom_{\cC}(F(Y),H(X,Y)).
$$
We still have to show that the category of algebras over $(\fun(\DD_1,\cC),\otimes_{\DD})$ is the category $\lax(\DD,\cC)$. Using again \ref{indmon}, it is enough to show that the full subcategory of algebras in $\wspan(\DD_1,\cC)$ spanned by the algebras whose underlying spans have the form $\DD_1 \overset{\Id{\DD_1}}{\longleftarrow} \DD_1 \overset{X}{\longrightarrow} \cC$ is equivalent to $\lax (\DD, \cC)$. Consider now a span $S=\DD_1 \overset{\Id{\DD_1}}{\longleftarrow} \DD_1 \overset{X}{\longrightarrow} \cC$ with a structure of algebra $(S,\mu_S,1_S)$. The multiplication $\mu_S \colon S \otimes S \to S$ is given by a diagram
$$
\begin{tikzcd}[row sep=small]
& \DD_1\times_{\DD_0} \DD_1 \arrow[rd,"\otimes_{\cC}\circ (S\times S)"{name=U}]  \arrow[ld,"\odot"'] \arrow[dd,"\odot"'] & \\
\DD_1 & & \cC \\
& \DD_1 \arrow[ur,"S"'] \arrow[Rightarrow,"\alpha"',shorten >=3mm,shorten <=3mm,from=U] \arrow[ul,"\Id{\DD_1}"] &
\end{tikzcd}
$$
and the unit $1_S \colon 1 \to S$ is given by a diagram 
$$
\begin{tikzcd}[row sep=small]
& \DD_0 \arrow[rd,"1_{\cC}"{name=U}]  \arrow[ld,"1_{\DD}"'] \arrow[dd,"1_{\DD}"'] & \\
\DD_1 & & \cC \\
& \DD_1 \arrow[ur,"S"'] \arrow[Rightarrow,"\beta"',shorten >=3mm,shorten <=3mm,from=U] \arrow[ul,"\Id{\DD_1}"] &
\end{tikzcd}
$$
and it is clear that the natural transformations $\alpha$ and $\beta$ endow $S\colon \DD_1 \to \cC$ with the structure of a lax functor from $\DD$ to $\cC$. If we consider now two such algebras $(S,\mu_{S},1_S)$ and $(T,\mu_T,1_T))$ in $\wspan(\DD_1,\cC)$, a morphism of algebras $(S,\mu_{S},1_S)\to (T,\mu_T,1_T))$ is given by a diagram
$$
\begin{tikzcd}[row sep=small]
& \DD_1 \arrow[rd,"S"{name=U}]  \arrow[ld,"\Id{\DD_1}"'] \arrow[dd,"\Id{\DD_1}"'] & \\
\DD_1 & & \cC \\
& \DD_1 \arrow[ur,"T"'] \arrow[Rightarrow,"\phi"',shorten >=3mm,shorten <=3mm,from=U] \arrow[ul,"\Id{\DD_1}"] &
\end{tikzcd}
$$
and the natural transformation $\phi\colon S \to T$ thus obtained produces a vertical natural transformation between the corresponding lax functors. This provides the equivalence between $\Algg(\fun(\DD_1,\cC),\otimes_{\DD})$ and $\lax (\DD ,\cC)$.
\end{proof}
\end{corollary}

\begin{remark}
With similar techniques, a convolution product on $\fun(\opo{\DD_1},\cC)$ can be constructed for other classes of double categories such as framed double categories. In this paper, we chose to concentrate only on target left fibrant double categories, as this class of double categories provides all the examples we are interested in.   
\end{remark}

\subsection{Composition product for \texorpdfstring{$\FF$}{F}-multicategories}
In this section, we use the convolution product for double categories that we constructed in the previous section to define a composition product for $\FF$-multicategories.
\begin{definition}\label{deftrf}
An $\FF$-multicategory $\MM$ is \emph{target right fibrant} if $t\colon \MM_1 \to \MM_0$ is a right fibration and if $\odot\colon \MM_1 \times_{\FF \MM_0}\FF \MM_1$ is an isofibration.
\end{definition}

\begin{proposition}
For every orbital pair $(\sO,\sT)$, the $\FF$-multicategory $\und{\Sigma_{(\sO,\sT)}}$ is a target right fibrant $\FF$-multicategory.
\begin{proof}
We prove first that the target functor $t\colon \Sigma_{(\sO,\sT)} \to \sO$ is a right fibration. If we consider any morphism $f\colon A \to B$ in $\sT$ such that $B$ belongs to $\sO$ and any morphism $g\colon C \to B$ in $\sO$, every pullback square 
$$
\begin{tikzcd}
D \arrow[r] \arrow[d] & \arrow[d,"f"] A \\
C \arrow[r,"g"'] & B
\end{tikzcd}
$$
provides a cartesian lift of $g\colon C \to t(f)$ in $\Sigma_{(\sO,\sT)}$. Since, by definition, morphisms in $\Sigma_{(\sO,\sT)}$ are pullback squares, this shows that $t$ is a right fibration.
We now prove that $\odot\colon \Sigma_{(\sO,\sT)}\times_{\FF \sO} \FF \Sigma_{(\sO,\sT)} \to \Sigma_{(\sO,\sT)}$ is an isofibration. If we consider any isomorphism 

\begin{equation}\tag{L}\label{doubsq}
\begin{tikzcd}
A \arrow[dd,"g"'] \arrow[r,"\simeq "] & B \arrow[d,"i"] \\
 & C\arrow[d,"j"] \\
D \arrow[r,"\simeq"'] & E,
\end{tikzcd}    
\end{equation}

in $\Sigma_{(\sO,\sT)}$, then the choice of any pullback square
$$
\begin{tikzcd}
F \arrow[r,"\simeq"] \arrow[d,"k"'] & C \arrow[d,"j"] \\
D \arrow[r,"\simeq "'] & E
\end{tikzcd}
$$ 
provides a lift of \ref{doubsq} given by the diagram
$$
\begin{tikzcd}
A \arrow[d,""'] \arrow[dd,bend right=30 ,"g"']\arrow[r,"\simeq "] & B \arrow[d,"i"] \\
F \arrow[r,"\simeq "] \arrow[d,"k"] & C\arrow[d,"j"] \\
D \arrow[r,"\simeq"'] & E,
\end{tikzcd} 
$$
seen as an isomorphism in $\Sigma_{(\sO,\sT)}\times_{\FF \sO} \FF \Sigma_{(\sO,\sT)}$.
\end{proof}
\end{proposition}

\begin{proposition}
The functor $\FF$ preserves cartesian fibrations, right fibrations, and isofibrations. In particular, if $\MM$ is a target right fibrant $\FF$-multicategory, then $\opo{\FF \MM}$ is a target left fibrant double category.
\begin{proof}
The first part of the statement is analogous to the proof of \cite[Lemma 2.20]{BonPer1}. Suppose now that $\MM$ is a target right fibrant $\FF$-multicategory.
The horizontal composition $\FF \MM_1 \times_{\FF \MM_0} \FF \MM_1 \to \FF \MM_1$ of $\FF \MM_1$ is given by $\FF \odot$ with $\odot\colon \MM_1 \times_{\FF \MM_0} \FF \MM_1 \to \MM_1$ the composition functor of $\MM$, using the isomorphism 
$$
\FF \MM_1 \times_{\FF \MM_0} \FF \MM_1 \simeq \FF \left( \MM_1 \times_{\FF \MM_0} \FF \MM_1\right).
$$
Moreover, the taget functor $\FF \MM_1 \to \FF \MM_0$ is the image under $\FF$ of the target functor $t\colon \MM_1 \to \MM_0$.
Finally, the result follows from the fact that the opposite functor of an isofibration is an isofibration and from the fact that the opposite functor of a right fibration is a left fibration. 
\end{proof}
\end{proposition}

\begin{proposition}\label{Fkan}
If the following diagram
$$
\begin{tikzcd}
 \opo{\cA} \arrow[r,"F"{name=U}]  \arrow[d,"\opo{G}"'] & \cV \\
\opo{\cB} \arrow[ur,"\operatorname{L}_GF"'] \arrow[Rightarrow,"\eta"',shorten >=3mm,shorten <=3mm,from=U] &
\end{tikzcd}
$$
exibits $\operatorname{L}_G F$ as the left Kan extension of $F$ along $\opo{G}$, then the following diagram
$$
\begin{tikzcd}
 \opo{\FF\cA }\arrow[r,"\opo{(\FF \opo{F})}"{name=U}]  \arrow[d,"\opo{\FF G}"'] 
 &\opo{(\FF \opo{\cV})} \arrow[r,"\widetilde{\Id{\cV}}"] & \cV  \\
\opo{\FF \cB} \arrow[rru,bend right=15,"\widetilde{\operatorname{L}_G F}"'] \arrow[ur,"\operatorname{L}_GF"'] \arrow[Rightarrow,"\FF \eta"',shorten >=3mm,shorten <=3mm,from=U] &
\end{tikzcd}
$$
is a left Kan extension diagram, with $\widetilde{(-)}$ as defined in \ref{opof}.
\begin{proof}
As suggested by \cite[Remark 2.25]{BonPer1}, the proof of this proposition is the same as the proof of \cite[Lemma 2.21]{BonPer1}. 
\end{proof}
\end{proposition}

\begin{corollary}\label{corocompprod}
For every target right fibrant $\FF$-multicategory, we obtain a monoidal structure $(\coll{\MM}(\cV),\circ_{\MM})$ on the category $\coll{\MM}(\cV)$ where $\circ_{\MM}$ is defined on two $\MM$-collections $S,T \colon \opo{\MM_1} \to \cV$ as the functor $S\circ_{\MM} T$ given by the left Kan extension
\begin{equation}\label{comp}\tag{A}
\begin{tikzcd}
\opo{\left(\MM_1 \times_{\FF \MM_0} \FF \MM_1\right)} \arrow[r,"S\times \widetilde{T}"] \arrow[d,"\opo{\odot_{\MM}}"'] & \cV \times \cV \arrow[ld,Rightarrow,shorten >=4mm,shorten <=4mm] \arrow[r,"\times"] & \cV.  \\
\opo{\MM_1} \arrow[rru,bend right=10,"S\circ_{\MM} T"'] & & 
\end{tikzcd}
\end{equation}
Moreover, the category of algebras in $(\coll{\MM}(\cV),\circ_{\MM})$ is $\oper{\MM}(\cV)$.
\begin{proof}
The category $\coll{\MM}(\cV)$ can be identified with the full subcategory of $\fun(\opo{\FF \MM_1},\cV)$ spanned by the functors $F\colon \opo{\FF \MM_1} \to \cV$ that preserve finite products. It is now enough to show that for all functors $S,T\colon \opo{\MM_1} \to \cV$, if we denote by $\widetilde{S},\widetilde{T}\colon \opo{\FF \MM_1} \to \cV$ the associated functors that preserve finite products, then the functor $\widetilde{S} \otimes_{\opo{\FF \MM}} \widetilde{T}$ is isomorphic to $\widetilde{S \circ_{\MM} T}$. The category $\FF\MM_1 \times_{\FF \MM_0} \FF \MM_1$ is equivalent to $\FF \left(\MM_1\times_{\FF \MM_0} \FF \MM_1 \right)$ and, through this equivalence, the functor $\opo{\left(\FF \MM_1 \times_{\FF \MM_0}\FF \MM_1  \right)} \to \opo{\FF \MM_1} \times \opo{\FF \MM_1} \overset{\widetilde{S}\times \widetilde{T}}{\longrightarrow} \cV \times \cV \overset{\times}{\longrightarrow} \cV$ is isomorphic to the composition of the image of $\opo{\left(\MM_1 \times_{\FF \MM_0} \FF \MM_1\right)} \overset{}{\longrightarrow} \opo{\MM_1} \times \opo{\FF \MM_1} \to \cV \times \cV \overset{\times}{\longrightarrow} \cV$ under $\opo{\FF(\opo{(-)})}$ with $\widetilde{\Id{\cV}} \colon \opo{\FF(\opo{\cV})} \to \cV$. The first part of the corollary now follows from \ref{Fkan} applied to the diagram \ref{comp}. The second part of the corollary follows from the fact that the inclusion $\coll{\MM}(\cV)\hookrightarrow \fun(\opo{\FF \MM_1},\cV)$ is monoidal and from \ref{algmon} using the description of $\MM$-operads in terms of lax functors.
\end{proof}
\end{corollary}

\begin{proposition}\label{trfcolour}
For every target right fibrant $\FF$-multicategory $\MM$ and for every right fibration $F\colon \cA \to \MM_0$, the $\FF$-multicategory $F^* \MM$ is also target right fibrant. In particular, for every $\MM$-set of colours $A$, the $\FF$-multicategory $\MM^A$ is target right fibrant if $\MM$ is target right fibrant.
\begin{proof}
We show first that $t\colon F^*\MM_1 \to \cA$ is a right fibration. The functor $t$ is by definition the projection $p_2\colon \MM_1 \times_{\left(\MM_0 \times \FF \MM_0\right)} \left( \cA \times \FF \cA \right) \overset{}{\longrightarrow} \cA$
that is a pullback of the functor $\MM_1 \times_{\FF \MM_0} \FF \cA \to \MM_0$.
The functor $\MM_1 \times_{\FF \MM_0} \FF \cA \to \MM_0$ is the composition
$$
\MM_1 \times_{\FF \MM_0} \FF \cA \to \MM_1 \overset{t_{\MM}}{\longrightarrow} \MM_0,
$$
where the first functor is a right fibration as a pullback of $\FF \cA \overset{\FF F}{\longrightarrow} \FF \MM_0$ and where $t_{\MM}$ is by assumptions a right fibration. It follows that $\MM_1 \times_{\FF \MM_0} \FF \cA \to \MM_0$ and $t$ are also right fibrations.
We now have to show that the functor 
$$
\odot_{F^*\MM}\colon \left(\MM_1 \times_{\FF \MM_0}\FF \MM_1\right)\times_{\left(\MM_0\times \FF \MM_0 \times \FF^2 \MM_0\right)}\left( \cA \times \FF \cA \times \FF^2 \cA \right) \to \MM_1 \times_{\left(\MM_0\times \FF \MM_0\right)} \left( \cA \times \FF \cA \right)
$$
is an isofibration. Note that we have an isomorphism 
$$
\MM_1 \times_{\left(\MM_0\times \FF \MM_0 \right)}\left(\cA \times \FF \cA\right)\simeq \MM_1 \times_{\left(\MM_0 \times \FF^2 \MM_0\right)} \left(\cA \times \FF^2 \cA \right)
$$
and that we have a pullback square
$$
\begin{tikzcd}
\left(\MM_1 \times_{\FF \MM_0}\FF \MM_1\right)\times_{\left(\MM_0\times \FF \MM_0 \times \FF^2 \MM_0\right)}\left( \cA \times \FF \cA \times \FF^2 \cA \right) \arrow[d,"G"'] \arrow[r] & \FF \cA \arrow[d,"\FF F"] \\
\left(\MM_1 \times_{\FF \MM_0}\FF \MM_1\right)\times_{\left(\MM_0\times  \FF^2 \MM_0\right)}\left( \cA \times  \FF^2 \cA \right) \arrow[r]  & \FF \MM_0.
\end{tikzcd}
$$
The functor $\odot_{F^* \MM}$ is the composition of $G$ with $\odot_{\MM}\times_{\Id{(\MM_0 \times \FF^2 \MM_0)}}(\cA \times \FF^2 \cA)$ and the result follows from the fact that $G$ is an isofibration as the pullback of an isofibration (using that $\FF$ preserves isofibrations) and from the fact that $\odot_{\MM}$ is an isofibration.
\end{proof}
\end{proposition}

\begin{corollary}\label{corocompcol}
For every target right fibrant $\FF$-multicategory $\MM$ and every $\MM$-set of colours $A$, the category of $A$-coloured $\MM$-collections in $\cV$ is endowed with the structure of a monoidal category for which the associated category of algebras is the category of $A$-coloured $\MM$-operads in $\cV$.
\begin{proof}
This follows directly from \ref{corocompprod} and \ref{trfcolour}.
\end{proof}
\end{corollary}

\subsection{Base change adjunctions}
In this section, we investigate when a multifunctor $\phi\colon \MM \to \NN$ between target right fibrant $\FF$-multicategories induces base change adjunctions between (coloured) $\MM$-operads and $\NN$-operads in $\cV$. In particular, we study base change adjunctions associated with changes of colours. 

\begin{proposition}\label{doublax}
Let $\phi\colon \DD \to \EE$ be a double functor between target left fibrant double categories. The functor $(\phi_1)^*\colon \fun (\EE_1,\cC) \to \fun (\DD_1,\cC)$ is naturally endowed with the structure of a lax functor. 
\begin{proof}
The following diagram 
$$
\begin{tikzcd}
\DD_1 \times_{\DD_0} \DD_1 \arrow[d,"\circ"'] \arrow[r,"\phi_1 \times_{\phi_0} \phi_1"] & \EE_1 \times_{\EE_0} \EE_1 \arrow[r,"S\times T"] \arrow[d,"\circ"'] & \arrow[ld,Rightarrow,shorten >=4mm,shorten <=4mm] \cC \times \cC \arrow[r,"\otimes_{\cC}"] & \cC \\
\DD_1 \arrow[r,"\phi_1"'] & \EE_1 \arrow[rru,bend right=10,"S\otimes_{\operatorname{D}}T"'] & &
\end{tikzcd}
$$
induces a morphism $\left(S\circ \phi_1\right)\otimes_{\operatorname{D}}\left(T\circ \phi_1 \right)\to  \left(S\otimes_{\operatorname{D}}T \right)\circ \phi_1$ by the universal property of left Kan extension and the following diagram
$$
\begin{tikzcd}
\DD_0 \arrow[r,"\phi_0"] \arrow[d,"1_{\DD}"'] & \EE_0 \arrow[d,"1_{\EE}"'] \arrow[r,"1_{\cC}"{name=U}] & \cC \\
\DD_1 \arrow[r,"\phi_1"'] & \EE_1 \arrow[Rightarrow,shorten >=3mm,shorten <=3mm,from=U]\arrow[ru,"1_{\fun(\DD_1,\cC)}"'] &
\end{tikzcd}
$$
induces a morphism $1_{\fun(\DD_1,\cC)} \to (\phi^*)(1_{\fun(\EE_1,\cC)})$. We thus obtain the structure of a lax functors on $(\phi_1)^*$.
\end{proof}
\end{proposition}

\begin{corollary}
Let $\phi \colon \MM \to \NN$ be a functor between target right fibrant $\FF$-multicategories. The functor $(\phi_1)^*\colon \coll{\NN}(\cV)\to \coll{\MM}(\cV)$ is naturally endowed with the structure of a lax functor. In particular, the functor $(\phi_1)_!\colon \coll{\MM}(\cV) \to \coll{\NN}(\cV)$ is naturally endowed with the structure of an oplax functor.
\begin{proof}
This follows directly from \ref{doublax} using that $(\phi_1)^*$ is the restriction of 
$$
\left(\opo{\FF \phi}_1\right)^*\colon \fun (\opo{\FF \NN_1},\cV) \to \fun (\opo{\FF\MM_1},\cC)
$$
to $\coll{\NN}(\cV)\subseteq \fun (\opo{\FF \NN_1},\cV)$ and $\coll{\MM}(\cV)\subseteq \fun (\opo{\FF \MM_1},\cV)$.
\end{proof}
\end{corollary}

\begin{proposition}\label{adjfib}
For every multifunctor $\phi\colon \MM \to \NN$, if for every $\NN$-set of colours $X$, the functor $\left(\phi^X\right)^*\colon \oper{\NN}^X(\cV)\to \oper{\MM}^{(\phi_0)^*X}(\cV)$ admits a right adjoint (respectively a left adjoint), then the functor $\phi^* \colon \operc{\NN}(\cV) \to \operc{\MM}(\cV)$ admits a right adjoint (respectively a left adjoint).
\begin{proof}
Suppose first that for every $\NN$-set of colours $X$, $\left(\phi^X\right)^*\colon \oper{\NN}^X(\cV) \to \oper{\MM}^{(\phi_0)^*X}(\cV)$ admits a right adjoint. If we consider the adjunction 
$$
(\phi_0)^* \colon \Psh(\NN_0) \rightleftarrows \Psh(\MM_0)\colon (\phi_0)_*,
$$
we have that $\Phi_{\phi}\colon \Phi_{\NN} \Rightarrow \Phi_{\MM}\circ (\phi_0)^*$ is a left morphism from $\Phi_\NN$ to $\Phi_{\MM}$ in the sense of \cite[Section 2.3]{grothmod} and we thus obtain the promised adjunction as described in \cite[Section 2.3]{grothmod}. Suppose now that for every $\NN$-set of colours $X$, $\left(\phi^X\right)^*\colon \oper{\NN}^X(\cV) \to \oper{\MM}^{(\phi_0)^*X}(\cV)$ admits a left adjoint. If we consider the adjunction 
$$
(\phi_0)_! \colon \Psh(\MM_0) \rightleftarrows \Psh(\NN_0)\colon (\phi_0)^*,
$$
then the natural transformation $\Phi_{\phi}\colon \Phi_{\NN} \Rightarrow \Phi_{\MM}\circ (\phi_0)^*$ produces a right morphism from $\Phi_{\MM}\circ (\phi_0)^*$ to $\Phi_{\phi}$ and we obtain the promised adjunction using \cite[Remark 2.3.5]{grothmod}.
\end{proof}
\end{proposition}

\begin{proposition}\label{critleft}
For every multifunctor $\phi\colon \MM \to \NN$ between target right fibrant $\FF$-multicategories such that $\phi_0$ is an isomorphism and such that $\phi_1$ is an isofibration, the functor 
$$
(\phi_1)_!\colon \left(\coll{\MM}(\cV) ,\circ_{\MM}\right) \to \left(\coll{\NN}(\cV),\circ_{\NN}\right)
$$
is strong monoidal.
\begin{proof}
Note first that because $\phi_0$ is an isomorphism, the morphism $(\phi_1)_!(1_{\coll{\MM}(\cV)})\to 1_{\coll{\MM}(\cV)}$ is an isomorphism. 
Consider now $S,T\colon \opo{\MM_1}\to \cV$ two $\MM$-collections in $\cV$. We have to show that the canonical natural transformation $(\phi_1)_!(S\circ_{\MM}T) \Rightarrow ((\phi_1)_!S)\circ_{\NN}((\phi_1)_!T)$ is an isomorphism. It is enough to show that the pasting of the following diagram
$$
\begin{tikzcd}[row sep=small]
\opo{\left(\MM_1 \times_{\FF \MM_0} \FF \MM_1\right)} \arrow[dd,""'] \arrow[r]  & \opo{\left(\MM_1\times \FF \MM_1 \right)}  \arrow[dd,""'] \arrow[rd,"S \times \widetilde{T} "{name=U}] & &   \\
 & & \cV \times \cV \arrow[r,"\times"] & \cV  \\
\opo{\left(\NN_1 \times_{\FF \NN_0} \FF \NN_1\right)} \arrow[dd,"\odot"']  \arrow[r] &  \opo{\left(\NN_1 \times \FF \NN_1 \right)} \arrow[ldd,Rightarrow,shorten >=4mm,shorten <=4mm]  \arrow[ur,""'] 
 \arrow[Rightarrow,""',shorten >=3mm,shorten <=3mm,from=U] & & \\
 &&& \\
 \opo{\NN_1} \arrow[uuurrr,bend right=20,""'] & & &
\end{tikzcd}
$$
is a left Kan extension. By definition the diagram 
$$
\begin{tikzcd}
\opo{\left(\NN_1 \times_{\FF \NN_0} \FF \NN_1\right)} \arrow[d,"\odot "'] \arrow[r] & \arrow[r,"S\times \widetilde{T}"] \opo{\left(\NN_1 \times \FF \NN_1\right)} \arrow[dl,Rightarrow,shorten >=4mm,shorten <=4mm] &  \cV \times \cV \arrow[r,"\times"] & \cV \\
\opo{\NN_1} \arrow[rrru,bend right=10,"S \circ_{\NN} T"'] &&&
\end{tikzcd}
$$
is a left Kan extension, and we have to check that the diagram 
$$
\begin{tikzcd}[row sep=small]
 \opo{\left(\MM_1 \times_{\FF \MM_0} \FF \MM_1\right)}\arrow[dd,""'] \arrow[r]  & \MM_1\times \FF \MM_1   \arrow[dd,""'] \arrow[rd,"S \times \widetilde{T} "{name=U}] & &   \\
 & & \cV \times \cV \arrow[r,"\times"] & \cV  \\
\opo{\left(\NN_1 \times_{\FF \NN_0} \FF \NN_1\right)} \arrow[r] &  \NN_1 \times \FF \NN_1   \arrow[ur,""'] 
 \arrow[Rightarrow,""',shorten >=3mm,shorten <=3mm,from=U] & & \\
\end{tikzcd}
$$
is also a left Kan extension diagram. Finally, using the pointwise formula for left Kan extension, it is enough to show that the functor 
$$    
\Phi\colon \opo{\left(\phi_1\times_{\phi_0}\FF\phi_1\right)}_{/(f,(f_1,\ldots,f_n))} \to \opo{\left(\phi_1 \times \FF \phi_1 \right)}_{/(f,(f_1,\ldots,f_n))}
$$
is final. Because $\phi_1$ is an isofibration, both the functors $\opo{\left(\phi_1 \times_{\phi_0} \FF \phi_1\right)}$ and $\opo{\left(\phi_1 \times \FF \phi_1\right)}$ are isofibrations and, moreover, $\opo{\left(\phi_1\times_{\phi_0} \FF \phi_1\right)}$ is a morphism of left fibrations as depicted by the diagram
$$
\begin{tikzcd}[column sep=tiny]
\opo{\left(\MM_1 \times_{\FF \MM_0}\FF\MM_1 \right)} \arrow[rr,"\opo{\left(\phi_1\times_{\phi_0} \FF \phi_1\right)}"] \arrow[rd,"t \circ p_1"'] & &\opo{\left(\NN_1 \times_{\FF \NN_0}\FF\NN_1 \right)} \arrow[ld,"t\circ p_1"] \\
 &\opo{\MM_0}&
\end{tikzcd}
$$
and the morphism $\opo{\left(\phi_1 \times \FF \phi_1\right)}$ is also a morphism of left fibrations as depicted by the following diagram
$$
\begin{tikzcd}[column sep=tiny]
\opo{\left(\MM_1 \times \FF \MM_1 \right)} \arrow[rr,"\opo{\left(\phi_1 \times \FF \phi_1\right)}"] \arrow[rd,"\opo{(t\times \FF t)}"'] & &\opo{\left(\NN_1 \times \FF\NN_1 \right)} \arrow[ld,"\opo{(t\times \FF t)}"] \\
 &\opo{\left(\MM_0 \times \FF \MM_0 \right)}.&
\end{tikzcd}
$$
It now follows from \ref{fiberfi} that the vertical functors of the square
$$
\begin{tikzcd}
\left(\opo{\left(\phi_1\times_{\phi_0}\FF\phi_1\right)}\right)^{-1}(f,(f_1,\ldots,f_n)) \arrow[r,"\simeq "] \arrow[d,hookrightarrow]& \left(\opo{\left(\phi_1\times\FF\phi_1\right)}\right)^{-1}(f,(f_1,\ldots,f_n)) \arrow[d,hookrightarrow] \\
\opo{\left(\phi_1\times_{\phi_0}\FF\phi_1\right)}_{/(f,(f_1,\ldots,f_n))} \arrow[r] & \opo{\left(\phi_1 \times \FF \phi_1 \right)}_{/(f,(f_1,\ldots,f_n))}
\end{tikzcd}
$$
are final and, because the upper horizontal functor is an equivalence, it follows that the bottom horizontal functor is also final.

\end{proof}
\end{proposition}

\begin{proposition}\label{summulticat}
The category $\multi{\FF}$ admits finite coproducts. Moreover, for all target right fibrant $\FF$-multicategories $\MM$ and $\NN$, $\MM \sqcup \NN$ is target right fibrant and the category $\oper{\MM\sqcup \NN}(\cV)$ is equivalent to the category $\oper{\MM}(\cV)\times \oper{\NN}(\cV)$.
\begin{proof}
Let $\MM$ and $\NN$ be two $\FF$-multicategories. We define $(\MM\sqcup\NN)_0$ and $(\MM \sqcup\NN)_1$ respectively by $\MM_0 \sqcup \NN_0$ and $\MM_1 \sqcup \NN_1$, the target functor is $t_{\MM}\sqcup t_{\NN}\colon \MM_1 \sqcup \NN_1 \to \MM_0 \sqcup \NN_0$, the identity functor is $1_{\MM}\sqcup 1_{\NN}\colon \MM_0 \sqcup \NN_0 \to \MM_1 \sqcup \NN_1$, the source functor is the composition 
$$
\MM_1 \sqcup \NN_1 \overset{s_{\MM}\sqcup s_{\NN}}{\longrightarrow} \FF(\MM_0) \sqcup \FF (\NN_0) \overset{(A,B)}{\longrightarrow} \FF(\MM_0\sqcup \NN_0)
$$
where the right hand side functor is the canonical functor given by the universal property of the coproduct, and finally, the composition is given by
$$
(\MM_1 \sqcup \NN_1)\times_{\FF\left(\MM_0 \sqcup \NN_0\right)}\FF (\MM_1 \sqcup \NN_1) \simeq \left( \MM_1 \times_{\FF \MM_0} \FF \MM_1\right)\sqcup \left(\NN_1 \times_{\FF \NN_0} \FF \NN_1\right)\overset{\odot_{\MM} \sqcup \odot_{\NN}}{\longrightarrow} \MM_1 \sqcup \NN_1
$$
where the first isomorphism follows from the universality of coproducts in $\Cat$ and from the fact that $\FF$ preserves pullbacks. It is clear that $\MM\sqcup\NN$ is target right fibrant using that both $\MM$ and $\NN$ are target right fibrant and that coproducts of right fibrations and isofibrations are respectively right fibrations and isofibrations. The last part of the statement follows from the fact that the equivalence $\coll{\MM\sqcup \NN}(\cV)\simeq \coll{\MM}(\cV)\times\coll{\NN}(\cV)$ is monoidal. 
\end{proof}
\end{proposition}

\begin{notation}
If $\MM$ is a target right fibrant $\FF$-multicategory and if $i\colon A \to B$ is a morphism of $\MM$-sets of colours, then we denote by $i^*$ the functor $(\MM^i)^*\colon \oper{\MM}^B(\cV) \to \oper{\MM}^A(\cV)$.
\end{notation}

\begin{construction}
Suppose that $\MM$ is a target right fibrant $\FF$-multicategory and let $i\colon A \hookrightarrow B$ be a monomorphism of $\MM$-sets of colours. If we write $i$ as an inclusion of summand $i\colon A \hookrightarrow A\sqcup C $, then the multifunctor $\MM^i\colon \MM^A \to \MM^B$ is the composition 
$$
\MM^A \overset{}{\hookrightarrow}\MM^A \sqcup \MM^C \overset{\left(\MM^{i_A},\MM^{i_C}\right)}{\longrightarrow} \MM^{A \sqcup C}
$$
where $i_A \colon A \hookrightarrow A\sqcup C$ and $i_C\colon C \hookrightarrow A \sqcup C$ are the canonical inclusions. The multifunctor 
$$
\left(\MM^{i_A},\MM^{i_C}\right)\colon \MM^A \sqcup \MM^C \to \MM^{A\sqcup C}
$$
satisfies the conditions of \ref{critleft} and it follows that the functor $i^* \colon \oper{\MM}^B(\cV) \to \oper{\MM}^A(\cV)$ admits a left adjoint $i_! \colon \oper{\MM}^A(\cV) \to \oper{\MM}^B (\cV)$ given by the composition
$$
\oper{\MM}^A(\cV) \overset{}{\longrightarrow} \oper{\MM}^A(\cV)\times \oper{\MM}^C(\cV) \overset{\left(\MM^{i_A}_1,\MM^{i_C}_1\right)_!}{\longrightarrow} \oper{\MM}^B(\cV)
$$
where the first functor is given component-wise by the identity of $\oper{\MM}^A(\cV)$ and by the functor constant on the initial object of $\oper{\MM}^C(\cV)$.
\end{construction}

\begin{corollary}\label{corex}
For every monomorphism $i\colon A \hookrightarrow B$ of $\MM$-sets of colours, the adjunction
$$
i_! \colon \oper{\MM}^A(\cV) \rightleftarrows \oper{\MM}^B (\cV) \colon i^*
$$
is coreflexive.
\begin{proof}
Using the previous construction, we can write
$$
i_! \colon \oper{\MM}^A(\cV) \rightleftarrows \oper{\MM}^B (\cV) \colon i^*
$$
as a composition of two adjunctions
$$
\begin{tikzcd}
\oper{\MM}^A(\cV) \arrow[r,shift left] & \arrow[l,shift left] \oper{\MM}^A(\cV)\times \oper{\MM}^C(\cV) \arrow[r,shift left,"(\MM^{i_A}_1\text{,}\MM_1^{i_C})_!"] & \oper{\MM}^B(\cV) \arrow[l,shift left,"(\MM_1^{i_A}\text{,}\MM_1^{i_C})^*"]
\end{tikzcd}
$$
and it is enough to show that both these adjunctions are coreflexive. It is clear that the left hand side adjunction is coreflexive. For the right hand side adjunction, this follows from the fact that the underlying monoidal adjunction
$$
\left(\MM_1^{i_A}\text{,}\MM_1^{i_C}\right)_! \colon \coll{\MM}^A(\cV)\times \coll{\MM}^C(\cV)  \rightleftarrows \coll{\MM}^B(\cV) \colon \left(\MM_1^{i_A}\text{,}\MM_1^{i_C}\right)^*
$$
is coreflexive as the functor $\left(\MM^{i_A}_1,\MM^{i_B}_1\right)$ is fully faithful. 
\end{proof}
\end{corollary}

\begin{proposition}\label{mate}
For every multifunctor $\phi\colon \MM \to \NN$ between target right fibrant $\FF$-multicategories such that $\phi_0$ is an isomorphism and for every monomorphism $i\colon A \hookrightarrow B$ of $\NN$-sets of colours, the mate

\begin{equation}\label{sqmate1}\tag{A}
\begin{tikzcd}
\oper{\NN}^B(\cV) \arrow[r,"(\phi^B)^*"]  & \oper{\MM}^B(\cV)  \\
\oper{\NN}^A(\cV) \arrow[u,"i_!"] \arrow[r,"(\phi^A)^*"'] & \oper{\MM}^A(\cV) \arrow[u,"i_!"'] \arrow[ul,Rightarrow,shorten >=2mm,shorten <=2mm]
\end{tikzcd}
\end{equation}

associated with the commutative square
\begin{equation}\label{sqmate2}\tag{B}
\begin{tikzcd}
\oper{\NN}^B(\cV) \arrow[d,"i^*"'] \arrow[r,"(\phi^B)^*"]  & \arrow[d,"i^*"] \oper{\MM}^B(\cV)  \\
\oper{\NN}^A(\cV)  \arrow[r,"(\phi^A)^*"'] & \oper{\MM}^A(\cV)  
\end{tikzcd}
\end{equation}

is an isomorphism.
\begin{proof}
If we write $i\colon A \hookrightarrow B$ as an inclusion of the form $i\colon A \hookrightarrow A \sqcup C$, the square \ref{sqmate2} can be decomposed into a composition of two squares as follows
\begin{equation}\tag{C}
\begin{tikzcd}
\oper{\NN}^{A\sqcup C}(\cV) \arrow[d,""'] \arrow[r,"(\phi^{A \sqcup C})^*"]  & \arrow[d,""] \oper{\MM}^{A\sqcup C}(\cV)  \\ 
\oper{\NN}^A(\cV) \times \oper{\NN}^C(\cV) \arrow[r] \arrow[d]  & \arrow[d] \oper{\MM}^A(\cV) \times \oper{\MM}^C(\cV)     \\
\oper{\NN}^A(\cV)  \arrow[r,"(\phi^A)^*"'] & \oper{\MM}^A(\cV)  
\end{tikzcd} 
\end{equation}
and the mate \ref{sqmate1} is the pasting of the diagram
\begin{equation}\tag{D}\label{sqmate4}
\begin{tikzcd}
\oper{\NN}^{A\sqcup C}(\cV)  \arrow[r,"(\phi^{A \sqcup C})^*"]  &  \oper{\MM}^{A\sqcup C}(\cV)  \\ 
\oper{\NN}^A(\cV) \times \oper{\NN}^C(\cV) \arrow[r] \arrow[u]  & \arrow[u] \oper{\MM}^A(\cV) \times \oper{\MM}^C(\cV)  \arrow[ul,Rightarrow,shorten >=5mm,shorten <=5mm]   \\
\oper{\NN}^A(\cV)  \arrow[r,"(\phi^A)^*"'] \arrow[u] & \oper{\MM}^A(\cV).  \arrow[u] \arrow[ul,Rightarrow,shorten >=5mm,shorten <=5mm]
\end{tikzcd} 
\end{equation}
It is now enough to show that the two mates of the diagram \ref{sqmate4} are isomorphisms. For the upper square, using that the vertical adjunctions are provided by monoidal adjunctions, it is enough to show that the following mate 
$$
\begin{tikzcd}
\fun\left(\opo{\left(\NN_1^{A\sqcup C}\right)},\cV\right)  \arrow[r,"(\phi_1^{A \sqcup C})^*"]  &   \fun\left(\opo{\left(\MM_1^{A\sqcup C}\right)},\cV\right) \\ 
\arrow[r] \arrow[u] \fun \left(\opo{\left(\NN_1^A \sqcup \NN_1^C\right)},\cV\right)  & \fun \left(\opo{\left(\MM_1^A \sqcup \MM_1^C\right)},\cV\right) \arrow[u]   \arrow[ul,Rightarrow,shorten >=5mm,shorten <=5mm]   
\end{tikzcd} 
$$
is an isomorphism. We have now to show that for every $(\mu,(b_1,\ldots,b_n),b_0)$ in $\MM_1^{A\sqcup C}$, the functor
$$
(\MM^{i_A}_1 , \MM^{i_C}_1)_{/(\mu,(b_1,\ldots,b_n),b_0)} \to (\NN^{i_A}_1 , \NN^{i_C}_1)_{/(\phi_1(\mu),(b_1,\ldots,b_n),b_0)}
$$
is final. Using that both $(\MM^{i_A}_1 , \MM^{i_C}_1)$ and $(\NN^{i_A}_1 , \NN^{i_C}_1)$ are cocartesian fibrations, it is enough to show that the upper horizontal functor of the square
$$
\begin{tikzcd}
(\MM^{i_A} , \MM^{i_C})^{-1}(\mu,(b_1,\ldots,b_n),b_0) \arrow[d,hookrightarrow] \arrow[r] & (\NN^{i_A} , \NN^{i_C})^{-1}(\phi_1(\mu),(b_1,\ldots,b_n),b_0) \arrow[d,hookrightarrow] \\
(\MM^{i_A} , \MM^{i_C})_{/(\mu,(b_1,\ldots,b_n),b_0)} \arrow[r] & (\NN^{i_A} , \NN^{i_C})_{/(\phi_1(\mu),(b_1,\ldots,b_n),b_0)}
\end{tikzcd}
$$
is an equivalence. This follows from the fact that the categories $(\MM^{i_A}_1 , \MM^{i_C}_1)^{-1}(\mu,(b_1,\ldots,b_n),b_0)$ and $(\NN^{i_A}_1 , \NN^{i_C}_1)^{-1}(\phi_1(\mu),(b_1,\ldots,b_n),b_0)$ are both empty if there exist $i$ and $j$ such that $b_i\in C$ and $b_j \in A$ and both a singleton if all $b_i$ belong to $A$ or if all $b_i$ belong to $C$.
\end{proof}
\end{proposition}

\subsection{A Monad for \texorpdfstring{$\MM$}{M}-operads} 
In what follows, $\MM$ is a target right fibrant $\FF$-multicategory. 
\begin{proposition}\label{leftclo}
The monoidal category $\left(\coll{\MM}(\cV),\circ_{\MM} \right)$ is left closed. 
\begin{proof}
The proof of this proposition is analogous to the proof that the convolution product is closed that we gave for \ref{algmon}.
Let $T$ be an $\MM$-collection in $\cV$. We have to show that the functor 
$$
(-)\circ_{\MM} T \colon \coll{\MM}(\cV) \to \coll{\MM}(\cV)
$$
admits a right adjoint. This functor can be obtained as the following composition
$$
\fun\left(\opo{\MM_1},\cV\right) \overset{}{\longrightarrow} \fun\left(\opo{\left(\MM_1 \times \FF \MM_1\right)},\cV\right) \overset{}{\longrightarrow} \fun\left(\opo{\left(\MM_1 \times_{\FF \MM_0} \FF \MM_1\right)},\cV\right) \overset{\odot_!}{\longrightarrow} \fun\left(\opo{\MM_1},\cV\right).
$$
The last functor is a left adjoint of $\odot^*$ and the second functor admits a right adjoint given by taking right Kan extensions. It remains to show that the first functor admits a right adjoint. As the category $\cV$ is cartesian closed, this functor is a left adjoint of the functor $\fun (\opo{\left(\MM_1 \times \FF \MM_1\right)},\cV) \to \fun (\opo{\MM_1},\cV)$ that sends a functor $H\colon \opo{\left(\MM_1 \times \FF \MM_1 \right)} \to \cV$ to the functor $\opo{\MM_1}\to \cV$ that sends an object $X$ of $\MM_1$ to the end
$$
\int_{Y \in \opo{\FF \MM_1}} \hom_{\cV}\left(\widetilde{T}(Y),H(X,Y)\right).
$$
\end{proof}
\end{proposition}

\begin{proposition}\label{closif}
The subcategory $\coll{\MM}(\cV)$ of $\fun(\opo{\FF \MM},\cV)$ is closed under sifted colimits. In particular, for every $\MM$-collection $T$ in $\cV$, the functor $\circ_{\MM}$ preserves sifted colimits in both variables. 
\begin{proof}
For every sifted category $I$ and every functor $F\colon I \to \coll{\MM}(\cV)$, the colimit of $F$ in the category $\fun (\opo{\FF \MM},\cV)$ is given by the composition
$$
\opo{\FF \MM} \to \fun (I,\cV) \overset{\colim}{\longrightarrow}  \cV.
$$
The first functor preserves products as $F$ factors by definition through the inclusion $\coll{\MM}(\cV)\hookrightarrow \fun (\opo{\FF \MM},\cV)$ and the second functor preserves finite products because $I$ is sifted and because $\cV$ is closed. The composition product $\circ_{\MM}$ is the restriction to $\coll{\MM}(\cV)$ of the convolution product $\otimes_{\MM}$ on $\fun (\opo{\FF \MM_1},\cV)$ and, as the convolution product $\otimes_{\MM}$ is closed, it follows that $\circ_{\MM}$ preserves sifted colimits in both variables.
\end{proof}
\end{proposition}

\begin{proposition}[\protect{\cite[Theorem 23.3]{Kel1}}]\label{Kelmon}
For every left closed cocomplete monoidal category $(\cC,\otimes)$ such that for all $X$ in $\cC$, the functor $X\otimes(-)\colon \cC \to \cC$ preserves filtered colimits, the category of algebras in $\cC$ is the category of algebras over a monad $\SSS$ on $\cC$. 
\end{proposition}

We recall the construction of the functor $\SSS \colon \cC \to \cC$.
We first define a functor $\SSS_n \colon \cC \to \cC$ by induction on $n$ as follows:
\begin{enumerate}
    \item the functor $\SSS_0$ is constant on the unit $I$ of $\cC$;
    \item for all $X$ in $\cC$, we define $\SSS_{n+1} X$ as $I \sqcup \left(X \otimes \SSS_{n}X \right)$;
\end{enumerate}
We construct now by induction a natural transformation $i_n \colon \SSS_n \to \SSS_{n+1}$ as follows:
\begin{enumerate}
    \item for all $X$ in $\cC$, the morphism $i_0$ is the canonical morphism $I \to I \sqcup X$;
    \item for all $X$ in $\cC$, the morphism $i_{n+1}\colon \SSS_{n+1} \to \SSS_{n+2}$ is the morphism
$$
\Id{I}\sqcup(\Id{X}\otimes i_n) \colon I \sqcup (X\otimes \SSS_n X) \to I \sqcup (X \otimes \SSS_{n+1}(X)).
$$
\end{enumerate}
Finally, for all $X$ in $\cC$, we define the image of $X$ under $\SSS$ as the colimit of the diagram
$$
\SSS_0 X \overset{i_0}{\longrightarrow} \SSS_1 X \overset{i_1}{\longrightarrow} \SSS_2 X\overset{i_2}{\longrightarrow} \cdots.
$$

\begin{corollary}
For every target right fibrant $\FF$-multicategory $\MM$, the category of $\MM$-operads in $\cV$ is the category of algebras over a monad $\SSS$ in $\coll{\MM}(\cV)$, that we denote by $\SSS_{\MM}$.
\begin{proof}
This follows directly from \ref{corocompprod}, \ref{leftclo}, \ref{closif}, and \ref{Kelmon}.  
\end{proof}
\end{corollary}

\begin{proposition}
The monad $\SSS_{\MM}$ preserves sifted colimits.
\begin{proof}
It is enough to show that for every non negative integer $n$, the functor $\SSS_{\MM}(n)$ preserves sifted colimits. The functor $\SSS_{\MM}(0)$ is constant and preserves sifted colimits as sifted categories are connected. Suppose now by induction that $\SSS_{\MM}(n)$ preserves sifted colimits. The functor $\SSS_{\MM}(n+1)$ is defined as the composition 
$$
\coll{\MM}(\cV) \overset{\Delta}{\longrightarrow} \coll{\MM}(\cV) \times \coll{\MM}(\cV) \overset{\Id{\cC} \times \SSS_{\MM}(n)}{\longrightarrow} \coll{\MM}(\cV) \times \coll{\MM}(\cV) \overset{\circ_{\MM}}{\longrightarrow} \coll{\MM}(\cV)
$$
and preserves sifted colimits as $\circ_{\MM}$ preserves sifted colimits in each variable and both $\Id{\cC}$ and $\SSS_{\MM}(n)$ preserve sifted colimits.
\end{proof}
\end{proposition}

\begin{corollary}\label{bicompop}
If $\cV$ is bicomplete, then the category $\oper{\MM}(\cV)$ is also bicomplete.    
\begin{proof}
The category $\fun (\opo{\MM_1},\cV)$ is bicomplete as the category $\cV$ is bicomplete and the result now follows from the fact that $\oper{\MM}(\cV)$ is the category of algebras over a monad on $\fun(\opo{\MM_1},\cV)$ that preserves sifted colimits. 
\end{proof}
\end{corollary}

\section{Homotopy theory of simplicial coloured \texorpdfstring{$\MM$}{M}-operads}\label{sec4}

In what follows, $\MM$ is a target right fibrant $\FF$-multicategory. A \emph{simplicial $\MM$-operad} (respectively a \emph{simplicial coloured $\MM$-operad}) is an $\MM$-operad (respectively a coloured $\MM$-operad) in $\sSet$. We will denote by $\soper{\MM}$ and $\soperc{\MM}$ the categories of simplicial $\MM$-operads and simplicial coloured $\MM$-operads respectively. We will also denote by $\oper{\MM}$ (respectively $\operc{\MM}$) the category of $\MM$-operads in $\Set$ (respectively coloured $\MM$-operads in $\Set$). In this section, our aim is to establish model structures on $\soper{\MM}$ and $\soperc{\MM}$, extending the classical model structures on simplicial (coloured) operads.

\subsection{Model structure on simplicial \texorpdfstring{$\MM$}{M}-operads}
In this subsection, we construct the  model structure on $\soper{\MM}$.

\begin{theorem}\label{THmodop}
There exists a model structure on the category of simplicial $\MM$-operads right induced by the free-forgetful adjunction
$$
F_{\MM}\colon \coll{\MM}(\sSet) \rightleftarrows s\oper{\MM}\colon U_{\MM}
$$
from the projective model structure on $\coll{\MM}(\sSet)$.
\begin{proof}
The proof is analogous to \cite[Theorem 3.2]{BM03}, and is using the transfer principle of \cite[Section 2.5]{BM03}. Note first that $U_{\MM}$ preserves filtered colimits because $\SSS_{\MM}$ does, and that $\soper{\MM}$ is bicomplete by \ref{bicompop}. By \ref{laxfunc}, the functor $\operatorname{Ex}^{\infty}$ induces a fibrant replacement functor on $s\oper{\MM}$ by post-composition. Finally, the interval $\Delta^0 \sqcup \Delta^0 \rightarrowtail \Delta^1 \overset{\sim}{\longrightarrow} \Delta^0$ induces a functorial path object $\cO \overset{\sim}{\longrightarrow} \cO^{\Delta^1}\twoheadrightarrow \cO \times \cO$ for fibrant operads in $s\oper{\MM}$ using \ref{laxfunc}.
\end{proof}
\end{theorem}

\begin{example}\label{excolormod}
For every $\MM$-set of colours $A$, using \ref{corocompcol}, we obtain a cofibrantly generated model structure on the category of simplicial $A$-coloured $\MM$-operads  $\soper{\MM}^A$.
\end{example}

\begin{example}\label{exgenmod}
For every orbital pair $(\sO,\sT)$, we thus obtain a model structure on the category $\soper{(\sO,\sT)}$ of simplicial $(\sO,\sT)$-operads. Note that in the case of the orbital pair $(\sO_G,\sO_{\cT})$ provided by a $G$-transfer system $\cT$ associated with a finite group $G$, this model structure is expected to correspond to the model structure constructed on \emph{partial genuine equivariant operads} in \cite[Theorem II]{BonPer1}. However, as we already pointed out in \ref{exgenop}, it is not known whether these two categories are equivalent.
\end{example}

\begin{example}
The $\FF$-multicategory $\und{\NN}$ is target right fibrant, and we thus obtain a model structure on the category of simplicial non-symmetric operads.  
\end{example}

\begin{proposition}\label{adjcolour}
For every multifunctor $\phi\colon \MM \to \NN$ between target fibrant $\FF$-multicategories, the functor $\phi^* \colon \soper{\NN} \to \soper{\MM}$ admits a left ajoint $\phi_!\colon \soper{\MM} \to \soper{\NN}$. In particular, the functor $\phi^* \colon \soperc{\NN} \to \soperc{\MM}$ also admits a left adjoint $\phi_! \colon \soperc{\MM} \to \soperc{\NN}$.
\begin{proof}
The category $\coll{\MM}(\sSet)$ is locally presentable as it is a category of presheaves. Moreover, the monad $\SSS_{\MM}$ preserves sifted colimits, and it follows that the category $\soper{\MM}$ is also locally presentable. Finally, the functor $\phi^*$ preserves limits and filtered colimits and the result now follows from the adjoint functor theorem \cite[Theorem 1.66]{AdRo}. The second part of the statement follows directly from \ref{adjfib}.
\end{proof}
\end{proposition}

\begin{proposition}
For every multifunctor $\phi\colon \MM \to \NN$ between target fibrant $\FF$-multicategories, the adjunction
$$
\phi_! \colon \soper{\MM} \rightleftarrows \soper{\NN}\colon \phi^*
$$
is a Quillen adjunction.
\begin{proof}
It is clear that $\phi^*$ is a right Quillen functor as weak equivalences and fibrations in simplicial $\MM$-operads are the morphisms which respectively forget to weak equivalences and fibrations of $\MM$-collections in $\sSet$. 
\end{proof}
\end{proposition}

\subsection{Model structure on simplicial coloured \texorpdfstring{$\MM$}{M}-operads}

In this subsection, we construct the model structure on $\soperc{\MM}$. The first step is to show that $\soperc{\MM}$ is bicomplete.

\begin{proposition}\label{simpob}
The category $\soperc{\MM}$ is equivalent to the subcategory of $\fun(\opo{\Delta},\operc{\MM})$ spanned by the simplical objects which are constant on colours.
\begin{proof}
Note first that we have an equivalence between $\soper{\MM}$ and $\fun(\opo{\Delta},\oper{\MM})$ for every $\FF$-multicategory $\MM$. This equivalence is functorial in $\MM$ and it follows that the category $\soperc{\MM}$ is equivalent to the Grothendieck construction of the functor $\opo{\Psh(\MM_0)} \to \CAT$ that sends $A$ to $\fun (\opo{\Delta},\oper{\MM}^A)$. Finally, it is clear that the Grothendieck construction of the later functor can be identified with the subcategory of $\fun(\opo{\Delta},\operc{\MM})$ spanned by these functors $\cO \colon \opo{\Delta} \to \operc{\MM}$ which are constant on colours.  
\end{proof}
\end{proposition}

We show now that the category $\soperc{\MM}$ is bicomplete.

\begin{definition}
The strict $2$-category $\Adj$ is the wide subcategory of $\CAT$ spanned by the functors $F\colon \cC \to \cD$ that admit a left adjoint. If $F\colon \opo{\cC} \to \Adj$ is a functor, we can in particular consider its Grothendieck construction $\Gamma_F \colon \int_{\cC}F \to \cC$.
\end{definition}

\begin{proposition}[\protect{\cite[Proposition 2.4.4]{grothmod}}]\label{coligro}
For every functor $F\colon \opo{\cC} \to \Adj$ such that $\cC$ is bicomplete and such that $F(X)$ is bicomplete for every $X$ in $\cC$, the category $\int_{\cC}F$ is also bicomplete. 
\end{proposition}

\begin{proposition}\label{bicomcol}
The category $\soperc{\MM}$ is bicomplete.
\begin{proof}
The result follows from \ref{coligro}, \ref{bicompop} and \ref{adjcolour} using that $\Psh(\MM_0)$ is bicomplete.
\end{proof}
\end{proposition}

If we denote by $\sCat$ the category of simplicial categories, then we have an adjunction
$$
i_! \colon \fun(\opo{\MM_0},\sCat) \rightleftarrows \soperc{\MM}\colon i^*
$$
provided by \ref{inducedcat} and \ref{adjcolour}
and we denote by $(i_X)_! \colon \sCat \rightleftarrows \soperc{\MM}\colon (i_X)^*$
the composition of adjunctions
$$
\begin{tikzcd}
\sCat \arrow[r,shift left=1.5,"\hom(-\text{,}X)\times(-)"]& \fun(\opo{\MM_0},\sCat) \arrow[r,shift left=1.5,"i_!"] \arrow[l,shift left=1.5,"\operatorname{ev}_X"] & \soperc{\MM}. \arrow[l,shift left=1.5,"i^*"]
\end{tikzcd}
$$

Recall that a functor $F\colon \cC \to \cD$ between simplicial categories is \emph{essentially surjective} (respectively an \emph{isofibration}) if $\pi_0F$ is essentially surjective (respectively if $\pi_0F$ is an isofibration).
We can now define the following classes of morphisms in $\soperc{\MM}$.  

\begin{definition}\label{weakfibop}
A morphism of simplicial coloured $\MM$-operads $(f,\phi)\colon (A,\cO) \to (B,\cP)$ is:
\begin{enumerate}
    \item \emph{fully-faithful} if $\phi \colon \cO \to f^* \cP$ is a weak equivalence in $\soper{\MM}^C$;
    \item \emph{essentially surjective} if for all $X$ in $\MM_0$, $(i_X)^*(f,\phi)$ is essentially surjective in $\sCat$;
    \item a \emph{weak equivalence} if it is fully-faithful and essentially surjective;
    \item a \emph{local fibration} if the morphism $\phi \colon \cO \to f^* \cP$ is a fibration in $\soper{\MM}^A$;
    \item an \emph{isofibration} if it is a local fibration and if for every $X$ in $\MM_0$, $(i_X)^*(f,\phi)$ is an isofibration in $\sCat$.
\end{enumerate}
A simplicial coloured $\MM$-operad $(A,\cO)$ is \emph{fibrant} if the unique morphism $(A,\sO) \to *$ is an isofibration, where $*$ is the final simplicial coloured $\MM$-operad. Equivalently, a simplicial coloured $\MM$-operad $(A,\sO)$ is fibrant if $\sO$ is fibrant in $\soper{\MM}^A$.
\end{definition}

We now want to show that the classes of weak equivalences and isofibrations are part of a model structure on $\soper{\MM}$. This model structure is cofibrantly generated and we introduce now its sets of generators.

\begin{notation}\label{notaclass}
Recall that $\FF \MM_0$ can be identified with the full subcategory of $\Psh(\MM_0)$ spanned by finite coproducts of representables.
For every $\MM$-set of colours $A$ that belongs to $\FF \MM_0$, we denote respectively by $I_A$ and $J_A$ the sets of generating cofibrations and generating trivial cofibrations of the model structure on $\oper{\MM}^A(\sSet)$ given by \ref{excolormod}.
Consider the following two classes of functors in $\sCat$:
\begin{enumerate}
    \item[$(C_2)$] the unique functor $\emptyset \to *$ where $\emptyset$ and $*$ are respectively the initial category and the final category;
    \item[$(A_2)$] the functors $* \to H$ which are fully faithful and essentially surjective, where $H$ is a simplicial category with set of objects $\{0,1\}$, with the property that $H_n$ is countable for each $n\ge 0$ and that $H$ is cofibrant in $\sCat_{\{0,1\}}$.
\end{enumerate}
 We now consider the following classes of morphisms in $\soperc{\MM}$:
\begin{enumerate}
    \item[$(J_1)$] the union of all the classes $J_A$ in $\soperc{\MM}$ for all $A$ in $\FF \MM_0$;
    \item [$(J_2)$] the union of the images of $A_2$ in $\soperc{\MM}$ under all the functors $(i_X)_!$ for all $X$ in $\MM_0$;
    \item[$(J)$] the union of $J_1$ and $J_2$;
    \item[$(I_1)$] the union of all the classes $I_A$ in $\soperc{\MM}$ for all $A$ in $\FF \MM_0$;
    \item[$(I_2)$] the union of all the images of $C_2$ in $\soperc{\MM}$ under all the functors $(i_X)_!$ for all $X$ in $\MM_0$;
    \item[$(I)$] the union of $I_1$ and $I_2$.
\end{enumerate}

\end{notation}

\begin{proposition}\label{wecompact}
A morphism $(f,\phi)\colon (A,\cO) \to (B,\cP)$ of simplicial coloured $\MM$-operads is a local weak equivalence (respectively a local fibration) if and only if for every morphism of $\MM$-sets of colours $g\colon C \to A$ with $C$ in $\FF \MM_0$, the morphism $g^*\phi\colon g^*\cO \to g^*f^*\cP$ is a weak equivalence (respectively a fibration) in $\soper{\MM}^C$.
\begin{proof}
We only show the result for local weak equivalences, the proof for local fibrations is analogous. Consider a morphism $(f,\phi)\colon (A,\cO)\to (B,\cP)$ of simplicial coloured $\MM$-operads. Suppose first that for every $g\colon C \to A$ with $C$ in $\FF \MM_0$, the morphism $g^*\phi$ is a weak equivalence. Let us show that $(f,\phi)$ is a local weak equivalence. We have to show that, for every $(\mu, (a_1,\ldots,a_n),a_0)$ in $\MM^A$, the map 
$$
\phi(\mu, (a_1,\ldots,a_n),a_0)\colon \cO(\mu,(a_1,\ldots,a_n),a_0) \to \cP(\mu,(f(a_1),\ldots,f(a_n)),f(a_0))
$$
is a weak equivalence in $\sSet$. Consider the object $C:=(a_0,\ldots,a_n)$ of $\FF \MM_0$ seen as an $\MM$-set of colours. 
We have a canonical morphism $g \colon C \to A$ using the Yoneda Lemma and the fact that every $a_i$ belongs to $A$. By assumptions $g^*(\phi)$ is a weak equivalence and it follows directly that $\phi(\mu, (a_1,\ldots,a_n),a_0)$ is also a weak equivalence. Suppose now that $(f,\phi)$ is a local weak equivalence. The morphism $\phi$ is a weak equivalence in $\soper{\MM}^A$ and it follows that for every morphism of $\MM$-sets of colours $g\colon C\to A$, the morphism $g^*(\phi)$ is a weak equivalence in $\soper{\MM}^C$, because $g^*$ preserves weak equivalences.
\end{proof}
\end{proposition}

\begin{proposition}\label{liftfib}
A morphism of simplicial coloured $\MM$-operads is a local fibration if and only if it has the right lifting property with respect to $J_1$. 
\begin{proof} 
Let $(f,\phi)\colon (A,\cO) \to (B,\cP)$ be a morphism of simplicial coloured $\MM$-operads. By definition of $J_1$, we have to show that $(f,\phi)$ is a local fibration if and only if it has the right lifting property with respect to $J_C$ for every $C$ in $\FF \MM_0$. Suppose first that $(f,\phi)$ has  the right lifting property with respect to $J_C$ for every $C$ in $\FF \MM_0$. Using \ref{wecompact}, we have to show that for every $C$ in $\FF \MM_0$ and every morphism $\iota\colon C\to A$, the morphism $\iota^*(\phi)$ is a fibration in $\soper{\MM}^C$. For every $g\colon \cQ \to \cR$ in $J_C$ consider any commutative square

\begin{equation}\label{lifsq}\tag{C}
\begin{tikzcd}
\cQ \arrow[r,"\psi"] \arrow[d,"g"'] & \iota^* \cO \arrow[d,"\iota^* \phi"] \\
\cR \arrow[r] & \iota^* f^* \cP
\end{tikzcd}
\end{equation}

in $\soper{\MM}^C$. If we now consider the associated square
$$
\begin{tikzcd}
(C,\cQ) \arrow[r,"(\iota \text{,}\psi)"] \arrow[d,"(\Id{C}\text{,}g)"'] & (A, \cO) \arrow[d,"(f\text{,}\phi)"] \\
(C,\cR) \arrow[r] & (B,  \cP)
\end{tikzcd}
$$
in $\soperc{\MM}$, we obtain the existence of a lifting map $(\iota,\eta)\colon (C,\cR) \to (A,\cO)$ and it follows that $\eta\colon \cR \to \iota^* \cO$ produces a lift for the square \ref{lifsq}. In particular, $\iota^*(\phi)$ has the right lifting property with respect to $J_C$ in $\soper{\MM}^C$ and is thus a fibration in $\soper{\MM}^C$. Suppose now that $(f,\phi)$ is a local fibration and consider a square
$$
\begin{tikzcd}
(C,\cQ) \arrow[r,"(\iota \text{,}\psi)"] \arrow[d,"(\Id{C}\text{,}g)"'] & (A, \cO) \arrow[d,"(f\text{,}\phi)"] \\
(C,\cR) \arrow[r] & (B,  \cP)
\end{tikzcd}
$$
with $g$ in $J_C$. Obtaining a lift for this square is equivalent as obtaining a lift for the square
\begin{equation}
\begin{tikzcd}
\cQ \arrow[r,"\psi"] \arrow[d] & \iota^* \cO \arrow[d,"\iota^* \phi"] \\
\cR \arrow[r] & \iota^* f^* \cP
\end{tikzcd}
\end{equation}
in $\soper{\MM}^C$. The result now follows from the fact that $\iota^* \phi$ is a fibration, as $\iota^*\colon \soper{\MM}^A \to \soper{\MM}^C$ is a right Quillen functor. 
\end{proof}
\end{proposition}

\begin{proposition}\label{lifttfib}
A morphism of simplicial coloured $\MM$-operads is a local trivial fibration if and only if it has the right lifting property with respect to $I_1$.
\begin{proof}
The proof is analogous to the proof of \ref{liftfib}.
\end{proof}
\end{proposition}

\begin{proposition}\label{isovsJ}
A morphism of simplicial coloured $\MM$-operads is an isofibration if and only if it has the right lifting property with respect to the set $J$.
\begin{proof}
A morphism of simplicial coloured $\MM$-operads $(f,\phi)\colon (A,\cO) \to (B,\cP)$ is an isofibration if and only if it is a local fibration and for every $X$ in $\MM_0$, $(i_X)^*(f,\phi)$ is an isofibration in $\sCat$. The result now follows from \ref{liftfib} and from the fact that isofibrations in $\sCat$ are precisely the functors with the right lifting property with respect to $A_2$ in $\sCat$, as proven in \cite{Juju}.
\end{proof}
\end{proposition}

\begin{proposition}\label{isowvsI}
A morphism of simplicial $\MM$-operads has the right lifting property with respect to $I$ if and only if it is an isofibration and a weak equivalence.
\begin{proof}
The proof of this proposition is analogous to the proof of \ref{isovsJ} and uses \ref{lifttfib} as well as the fact that essentially surjective isofibrations in $\sCat$ are precisely the functors with the right lifting property with respect to $A_2$ in $\sCat$, as proven in \cite{Juju}.
\end{proof}
\end{proposition}

Our goal is now to prove that the domains of morphisms in $J$ and $I$ are small relative to $J$-cells and $I$-cells respectively.

\begin{proposition}\label{filteredmono}
For every multifunctor $\phi \colon \MM \to \NN$ such that $\phi_0\colon \MM_0 \to \NN_0$ is an isomorphism, the functor $\phi^* \colon \soperc{\NN} \to \soperc{\MM}$ preserves filtered colimits of diagrams $F\colon I \to \soperc{\NN}$ such that, for every $f\colon i \to j$ in $I$, the morphism $F(f) \colon F(i) \to F(j)$ is a monomorphism on colours.
\begin{proof}
If we denote by $(A_i,\cO_i)$ the image of $i$ under $F$, the colimit of $F$ is the simplicial coloured $\NN$-operad $(A,\cO)$ with $A=\colim_{i\in I}A_i$ and $\cO=\colim_i (\alpha_i)_! \cO_i$ where $\alpha_i\colon A_i \to A$ is the coprojection. We now have that 
$$
\phi^* (A,\cO)=\left(A, \left(\phi^A\right)^*\left(\colim_i (\alpha_i)_! \cO_i\right)\right)
$$
and $\left(\phi^A\right)^*(\colim_i (\alpha_i)_! \cO_i)\simeq \colim_i\left(\left(\phi^A\right)^*(\alpha_i)_!\cO_i\right)$ because $\left(\phi^A\right)^*$ preserves filtered colimits. Finally, using \ref{mate}, we have that the canonical morphism  $\left(\phi^X\right)^*(\alpha_i)_!\to (\alpha_i)_!\left(\phi^X\right)^*$ is an isomorphism and the result follows.
\end{proof}
\end{proposition}

\begin{remark}\label{inclucolours}
Note that the morphisms in $I$ and $J$ are monomorphisms on colours. Because pushout of morphisms of simplicial coloured $\MM$-operads that are monomorphisms on colours are also monomorphisms on colours, it follows that $I$-cells and $J$-cells are sequential colimits of morphisms that are monomorphisms on colours. It follows that \ref{filteredmono} applies to $I$-cells and $J$-cells.
\end{remark}

\begin{proposition}\label{smalladj}
Let $F\colon \opo{\cC} \to \operatorname{AdjCat}$ be a functor such that $\cC$ and $F(X)$ admit sequential colimits for every $X$ in $\cC$, and such that, for every morphism $f\colon X \to Y$ in $\cC$, the functor $f^*\colon F(Y) \to F(X)$ preserves sequential colimits. If $(Y,y)$ is an object of $\int_{\cC}F$ such that $Y$ is small in $\cC$ and such that $y$ is small in $F(Y)$, then $(Y,y)$ is small relative to the sequential colimits
\begin{equation}\tag{A}\label{seqcoco}
(X_0,x_0) \overset{(f_0,\phi_0)}{\longrightarrow} (X_1,x_1) \overset{(f_1,\phi_1)}{\longrightarrow} \cdots
\end{equation}
such that the adjunction 
$$
(\alpha_i)_! \colon F(X_i) \rightleftarrows F(\colim_i X_i)\colon (\alpha_i)^*
$$ 
associated with the coprojection $\alpha_i \colon X_i \to \colim_{i} X_i$ is coreflexive.
\begin{proof}
We have to show that every morphism $(g,\psi)\colon (Y,y) \to \colim_{i} (X_i,x_i)$ factors through a morphism $(g',\psi')\colon (Y,y) \to (X_k,x_k)$ for some $k$. If we denote by $(X,x)$ the colimit of \ref{seqcoco}, then, using \cite[Proposition 2.4.4]{grothmod}, we have that 
$$
(X,x)=(\colim_{i}X_i, \colim_i (\alpha_i)_! x_i )
$$
with $\alpha_i\colon X_i \to \colim_i X_i$ the coprojection. By assumptions, $Y$ is small in $\cC$ and it follows that $g \colon Y \to \colim_i X_i$ factors through a morphism $g''\colon Y \to X_m$ for some $m$.
The object $y$ is small in $F(Y)$, and using that $g^*$ preserves sequential colimits, the morphism $\psi\colon y \to \colim g^* (\alpha_i)_! x_i$ factors through a morphism $\psi'' \colon y \to g^* (\alpha_l)_! x_l$ for some $l$. If we denote by $k$ the maximum of $m$ and $l$, then we obtain that $g$ factors through a morphism $g'\colon Y \to X_k$ and that $\psi$ factors through a morphism $\psi'\colon y \to g^*(\alpha_k)_!x_k$. Moreover, as we have that $g=\alpha_k \circ g'$, we obtain that $\psi$ actually factors through a morphism $\psi'\colon y \to (g')^* (\alpha_k)^* (\alpha_k)_!x_k\simeq (g')^* x_k$ using that $((\alpha_k)_!, (\alpha_k)^*)$ is coreflexive by assumptions. It now follows that $(g,\psi)$ factors through a morphism $(g',\psi')\colon (y,Y) \to (X_k,x_k)$.

\end{proof}
\end{proposition}

\begin{proposition}\label{small}
The domain of morphisms in $I$ and $J$ are small relative to $I$-cells and $J$-cells respectively.
\begin{proof}
The domains of the morphisms in $C_2$ and $A_2$ are small in $\sCat$ and, because the functor $(i_X)^*\colon \soperc{\MM}\to \fun (\opo{\MM_0},\sCat)$ preserves sequential colimits of morphisms that are monomorphisms on colours, the functor $(i_X)_!$ preserves small objects relative to this class of sequential colimits. Using \ref{inclucolours}, the domain of the morphisms in $I_2$ and $J_2$ in $\soperc{\MM}$ are small relative to $I$-cells and $J$-cells respectively.
We now show that, for every $A$ in $\FF \MM_0$, the domains of morphisms in $I_A$ and $J_A$ are small relative to $I$-cells and $J$-cells. By assumption, $A$ is small in $\Psh(\MM_0)$ and we already know that the domains of maps in $I_A$ and $J_A$ are small in $\soper{\MM}^A$. The result now follows from \ref{inclucolours}, \ref{filteredmono}, and \ref{corex}.
\end{proof}
\end{proposition}

\begin{proposition}\label{fibequi}
Let $(A,\cO)$ be a fibrant simplicial coloured $\MM$-operad and consider $(\mu, (a_1,\ldots,a_n),a_0)$ and $(\mu,(b_1,\ldots,b_n),b_0)$ two objects of $\MM^A$ associated to the same multimorphism $\mu$ in $\MM$ with source $(X_1,\ldots,X_n)$ and target $X_0$. Suppose $a_i \sim b_i$ in $(i_{X_i})^*(A,\cO)$ for all $i=0,\ldots,n$. The Kan complexes $\cO(\mu,(a_1,\ldots,a_n),a_0)$ and $\cO(\mu,(b_1,\ldots,b_n),b_0)$ are homotopy equivalent. Moreover, any morphism of coloured $\MM$-operads $(f,\phi)\colon (A,\cO) \to (B,\cP)$ induces functorially a homotopy equivalence $\cP(\mu,(f(a_1),\ldots,f(a_n)),f(a_0)) \to \cP(\mu,(f(b_1),\ldots ,f(b_n)),f(b_0))$.
\begin{proof}
By assumptions, we have maps of simplicial sets $\alpha_i \colon \Delta^0 \to \cO(1_{X_i},b_i,a_i)$
and $\beta_i \colon \Delta^0\to \cO(1_{X_i},a_i,b_i)$ for all $i=1,\ldots,n$ that exhibit the equivalences $a_i \sim b_i$ in $(i_{X_i})^*(A,\cO)$.
We thus obtain a map of simplicial sets  $\theta\colon  \cO(\mu,(a_1,\ldots,a_n),a_0)\to \cO(\mu,(b_1,\ldots,b_n),b_0)$ given by the composition
$$
\cO(\mu,(a_1,\ldots,a_n),a_0) \overset{\alpha'}{\longrightarrow} \prod_{i=1}^n\cO(1_{X_i},b_i,a_i)\times  \cO(\mu,(a_1,\ldots,a_n),a_0)\times \cO(1_{X_0},a_0,b_0) \overset{\circ}{\longrightarrow} \cO(\mu,(b_1,\ldots,b_n),b_0) 
$$
with $\alpha'=((\alpha_i),\Id{\cO(\mu,(a_1,\ldots,a_n),a_0)},\beta_0)$ as well as a map $\eta \colon  \cO(\mu,(b_1,\ldots,b_n),b_0)\to \cO(\mu,(a_1,\ldots,a_n),a_0)$ given by the composition
$$
\cO(\mu,(b_1,\ldots,b_n),b) \overset{\beta'}{\longrightarrow} \prod_{i=1}^n \cO(1_{X_i},a_i,b_i)\times \cO(\mu,(b_1,\ldots,b_n),b_0)\times \cO(1_{X_0},b_0,a_0) \overset{\circ}{\longrightarrow} \cO(\mu,(a_1,\ldots,a_n),a_0)
$$
with $\beta'=((\beta_i),\Id{},\alpha_0)$. Using that $\alpha_0$ and $\beta_0$ and $\alpha_i$ and $\beta_i$ form equivalences, we deduce that $\theta$ and $\eta$ form the desired homotopy equivalence. The second part of the statement follows from the fact that the images of $\alpha_i$ and $\beta_i$ under $(f,\phi)$ provide by the same construction the promised homotopy equivalence.
\end{proof}
\end{proposition}

\begin{proposition}\label{23}
The class of weak equivalences in $\soperc{\MM}$ satisfies the two-out-of-three property.
\begin{proof}
Using that the functor $\operatorname{Ex}^{\infty}$ induces a fibrant replacement functor on $\soperc{\MM}$, it is enough to show the result for fibrant simplicial coloured $\MM$-operads. Consider now two morphisms of coloured simplicial $\MM$-operads
$$
(A,\cO) \overset{(f,\phi)}{\longrightarrow} (B,\cP) \overset{(g,\psi)}{\longrightarrow} (C,\cQ).
$$
The two-out-of-three property is true for essentially surjective functors as it is already proven in the category $\sCat$ and we must now show that if two out of three of the morphisms $(f,\phi)$, $(g,\psi)$ and $(g,\psi)\circ (f,\phi)$ are weak equivalences, then the third is a local weak equivalence. If $(f,\phi)$ and $(g,\psi)$ are both weak equivalences, then $f^* (\psi)\circ \phi$ is a local weak equivalence in $\soper{\MM}^A$ because $f^*$ preserves weak equivalences. Similarly, if the composition $(g,\psi)\circ (f,\psi)$ is a weak equivalence and if $(g,\psi)$ is a weak equivalence, then using the two-out-of-three for weak equivalences in $\soper{\MM}^A$, we deduce that $(f,\psi)$ is a local weak equivalence. Finally, we have to show that if both the composition $(g,\psi)\circ (f,\phi)$ and $(f,\phi)$ are weak equivalences, then $(g,\psi)$ is a local weak equivalence. Using again the two-out-three of weak equivalences in $\soper{\MM}^A$, we deduce that $f^*(\psi)$ is a weak equivalence in $\soper{\MM}^A$ and we want to show that $\psi$ is a weak equivalence in $\soper{\MM}^B$. For every $(\mu,(b_1,\ldots,b_n),b_0)$ in $\MM^B$, we want to show that 
$\psi(\mu,(b_1,\ldots,b_n),b_0) \colon  \cP(\mu,(b_1,\ldots,b_n),b_0) \to \cQ(\mu,(g(b_1),\ldots,g(b_n)),g(b_0))$ is a weak equivalence in $\soper{\MM}^B$. Because $i^*(f,\phi)$ is essentially surjective, we can find colours $a_0,\ldots,a_n$ in $A$ such that $b_i\sim f(a_i)$ for every $i=0,\ldots,n$. Finally, using \ref{fibequi}, we obtain a commutative square
$$
\begin{tikzcd}
\cP(\mu,(b_1,\ldots,b_n),b_0) \arrow[r,"\simeq"] \arrow[d,"\psi(\mu\text{,}(b_1\text{,}\ldots\text{,}b_n)\text{,}b_0)"'] & \cP(\mu,(f(a_1),\ldots,f(a_n)),f(a_0)) \arrow[d,"\simeq"] \\
\cQ(\mu,(g(b_1),\ldots,g(b_n)),g(b_0)) \arrow[r,"\simeq"'] & \cQ (\mu,(gf(a_1),\ldots,gf(a_n)),gf(a_0) )
\end{tikzcd}
$$
where the horizontal arrows are homotopy equivalences and where the right vertical arrow is a weak equivalence because $f^*\psi$ is a weak equivalence. Finally, we deduce that $\psi$ is a weak equivalence, and, in particular, that $(g,\psi)$ is a local weak equivalence.
\end{proof}
\end{proposition}

\begin{proposition}\label{retract}
The class of weak equivalences in $\soperc{\MM}$ is closed under retracts. 
\begin{proof}
We show that both local weak equivalences and essentially surjective morphisms are closed under retracts. For essentially surjective morphisms, this follows from the fact that this property in verified in the category $\sCat$. Suppose now that a morphism $(f,\phi)\colon (A,\cO) \to (B,\cP)$ is a retract of a local weak equivalence $(g,\psi)\colon (C,\cQ) \to (D,\cR)$. If we consider the associated diagram
$$
\begin{tikzcd}
(A,\cO) \arrow[r] \arrow[d,"(f\text{,}\phi)"'] & (C,\cQ) \arrow[d,"(g\text{,}\psi)"] \arrow[r] & (A,\cO)\arrow[d,"(f\text{,}\phi)"] \\
(B,\cP) \arrow[r] & (D,\cR) \arrow[r,"(h\text{,}\theta)"'] & (B,\cP),
\end{tikzcd}
$$
we deduce that $\phi$ is a retract of $h^*\psi$ and the result follows from the fact that $h^*\psi$ is a weak equivalence in $\soper{\MM}^B$ and from the fact that weak equivalences in $\soper{\MM}^B$ are closed under retracts. 
\end{proof}
\end{proposition}

\begin{proposition}\label{filcol}
For every sequence 
\begin{equation}\tag{B}\label{seqo2}
(A_0,\cO_0) \overset{(f_0,\phi_0)}{\longrightarrow} (A_1,\cO_1) \overset{(f_1,\phi_1)} {\longrightarrow} \cdots  
\end{equation}
of morphisms in $\soperc{\MM}$ that are monomorphisms on colours and weak equivalences, the canonical morphism $\alpha_0 \colon  (A_0,\cO_0) \to \colim_i (A_i,\cO_i)$ is a weak equivalence and a monomorphism on colours.
\begin{proof}
A sequential colimit of monomorphisms in $\Psh(\MM_0)$ is a monomorphism, and the second part of the statement follows.
We now show that the canonical morphism $\alpha_0\colon (A_0,\cO_0) \to \colim_i (A_i,\cO_i)$ is a weak equivalence. Using \ref{filteredmono} on the functor 
$$
i^* \colon \soperc{\MM} \to \fun(\opo{\MM_0},\sCat)
$$
and using the fact that weak equivalences in $\sCat$ are closed under sequential colimits, we already know that $\alpha_0 \colon (A_0,\cO_0) \to \colim_{i}(A_i,\cO_i)$ is essentially surjective, and we must prove that it is also fully faithful. The colimit of the sequence \ref{seqo2} is the simplicial $\MM$-operad 
$$
(\colim_i A_i,\colim_i (\alpha_i)_! \cO_i )
$$
and we have to prove that the canonical morphism $\phi_0 \colon \cO_0 \to  (\alpha_0)^*  (\colim_i (\alpha_i)_! \cO_i )$ is a weak equivalence in $\soper{\MM}^{A_0}$. If we denote by $\alpha_{i,j}\colon A_i \to A_j$ the morphism provided by the sequence $\ref{seqo2}$ whenever $j>i$, we have that $\alpha_j\circ \alpha_{i,j}= \alpha_i$ and using that $((\alpha_i)_!,(\alpha_i)^*)$ is coreflexive, we obtain that 
$$
(\alpha_0)^*  (\colim_i (\alpha_i)_! \cO_i ) \simeq \colim_i (\alpha_0)^*(\alpha_i)_! \cO_i\simeq  \colim_i \alpha_{0,i}^* \cO_i.
$$
Finally, the morphism $\phi$ is the canonical morphism $\cO_0 \to \colim_i \alpha_{0,i}^* \cO_i$ associated with the sequential colimit
$$
\cO_0=\colim_i \alpha_{0,0}^* \cO_0\overset{}{\longrightarrow} \colim_i \alpha_{0,1}^* \cO_1 \overset{}{\longrightarrow} \cdots
$$
in $\soper{\MM}^{A_0}$. This morphism is a weak equivalence, as it is a sequential colimit of weak equivalences in $\soper{\MM}^{A_0}$.
\end{proof}
\end{proposition}

\begin{proposition}\label{pushj}
Let $u\colon K \to H$ be a fully faithful functor between small categories such that the set of objects of $H$ is $\{0,1\}$, while the only object of $K$ is $0$ (such that $u(0)=0$). For every presheaf $F\colon \opo{\MM_0} \to \Set$ and every pushout square
$$
\begin{tikzcd}
i_!\left(F \times K \right) \arrow[r,"v"] \arrow[d,"i_!(\Id{F}\times u)"'] & (A,\cO) \arrow[d,"(f\text{,}\phi)"] \\
i_!\left(F \times H \right) \arrow[r] & (B,\cP)
\end{tikzcd}
$$
in the category $\operc{\MM}$, the morphism $(f,\phi)$ is fully faithful. If moreover $u$ is an equivalence, then $(f,\phi)$ is also an equivalence.
\begin{proof}
The proof of this proposition is similar to the proof of \cite[Lemma 1.29]{cismo}.
Note first that we see here both $H$ and $K$ as the objects of $\fun(\opo{\MM_0},\sCat)$ given by the associated constant functors.
The proof of this proposition relies on an explicit description of the $\MM$-operad $(B,\cP)$. We denote by $w$ the parameterized functor $F\times K \to i^*(A,\cO)$ associated to $u$ by adjunction and, for every $X$ in $\MM_0$ and $x$ in $F(X)$, we denote by $w_x\colon K \to i^*(A,\cO)(X)$ the associated functor.
For every $X$ in $\MM_0$ and every $x$ in $F(X)$, we denote by $x_0$ the image in $A(X)$ of $x$ under $w_x$.
The $\MM$-set of colours $B$ is the $\MM$-set $A \sqcup F$ and for all $X$ in $\MM_0$ and $x$ in $F(X)$, we denote by $x_1$ the corresponding element of $B(X)$ through the inclusion $F(X) \hookrightarrow A(X) \sqcup F(X)= B(X)$. Consider now a multimorphism $\mu$ of $\MM$ whose source and target are respectively $(X_1,\ldots,X_n)$ and $X_0$. For every $(\mu,(b_1,\ldots,b_n),b_0)$ in $\MM^B$, we define
$$
c_i=\begin{cases}
b_i & \text{if}~~ b_i \in  A(X_i); \\
x_0 & \text{if}~~b_i=x\in F(X_i).
\end{cases}
$$
We now define $\cP(\mu, (b_1,\ldots,b_n),b_0)$ as follows:
\begin{enumerate}
\item Suppose that $c_0$ belongs to $A(X_0)$. We denote by $\cQ(\mu, (b_1,\ldots,b_n),b_0)$ the set
$$
\cO (\mu, (c_1,\ldots,c_n),c_0)\times \prod_{i, b_i \in F(X_i),i\neq 0 } H(1,0).  
$$
We can consider on $\cQ(\mu, (b_1,\ldots,b_n),b_0)$ the equivalence relation generated by
$$
(f,(h_1\circ g_1,\ldots,h_m \circ g_m))\sim (f\circ (l_1,\ldots,l_n),(g_1,\ldots,g_m)) 
$$
where $(f,(g_1,\ldots,g_m))$ belongs to $\cQ(\mu, (b_1,\ldots,b_n),b_0)$, where $h_i \in K(0,0)$ for $i=1,\ldots,m$ with $m=|\{b_i~|~b_i \in A(X_i)\}|$ and where $l_i$ is defined by
$$
l_i=\begin{cases}
1_{b_i} & \text{if}~~ b_i \in  A(X_i); \\
w_{b_i}(h_i)  & \text{if}~~b_i\in F(X_i).
\end{cases}
$$
We define the set $\cP(\mu, (b_1,\ldots,b_n),b_0)$ as the quotient of $\cQ(\mu, (b_1,\ldots,b_n),b_0)$ by this equivalence relation;
\item We suppose now that $c_0$ belongs to $F(X_0)$. We denote by $\cQ(\mu, (b_1,\ldots,b_n),b_0)$ the set
$$
H(0,1)\times \cO (\mu, (c_1,\ldots,c_n),c_0)\times \prod_{i, b_i \in F(X_i),i\neq 0 } H(1,0).  
$$
We can consider on $\cQ(\mu, (b_1,\ldots,b_n),b_0)$ the equivalence relation generated by
$$
(g_0\circ h_0,f,(h_1\circ g_1,\ldots,h_m \circ g_m))\sim (g_0,l_0 \circ f\circ (l_1,\ldots,l_n),(g_1,\ldots,g_m)) 
$$
where $(g_0,f,(g_1,\ldots,g_m))$ belongs to $\cQ(\mu, (b_1,\ldots,b_n),b_0)$, where $h_i \in K(0,0)$ for $i=1,\ldots,m$ and where $h_0 \in H(0,1)$ with $m=|\{b_i~|~b_i \in A(X_i)\}|$ and where $l_i$ is defined by
$$
l_i=\begin{cases}
1_{b_i} & \text{if}~~ b_i \in  A(X_i) ; \\
w_{b_i}(h_i)  & \text{if}~~b_i \in F(X_i); \\

\end{cases}
$$
We define the set $\cP(\mu, (b_1,\ldots,b_n),b_0)$ as the quotient of $\cQ(\mu, (b_1,\ldots,b_n),b_0)$ by this equivalence relation.
\end{enumerate}
We describe the structure of $B$-coloured $\MM$-operad of $\cP$. We only explain how to compose operations corresponding to case (2), the other possible compositions can be defined in a similar way. Consider composable operations $(\lambda_0,(\lambda_1,\ldots,\lambda_n ))$ in $\MM^B_1 \times_{\FF \MM^B_0} \FF \MM^B_1$ such that the source of $\lambda_0$ is the sequence $(b_1,\ldots,b_n)$ of elements of $B$. For all $i=0,\ldots,n$, let $p_i$ be an element of $\cP(\lambda_i)$ represented by $(g_i^0,f_i,g_i^1,\ldots,g_i^{m_i})$. We define $p_0\circ (p_1,\ldots,p_n)$ as the operation of $\cP$ that is represented by
$$
(g_0^0, f_0 \circ (k_1 \circ l_1,\ldots,k_n\circ l_n) \circ (f_1,\ldots,f_n)  ,(g_1^1,\ldots,g_1^{m_1},g_2^1,\ldots, g_n^{m_n}))
$$
where $l_i$ and $k_i$ are respectively defined for $i=1,\ldots,n$ by
$$
l_i=\begin{cases}
g_0^i  & \text{if} ~~ b_i \in F(X_i) \\
1_{b_i} & \text{if} ~~ b_i \in A(X_i)
\end{cases}
~~\text{and}~~
k_i=\begin{cases}
g_i^0 & \text{if} ~~ b_i \in F(X_i) \\
1_{b_i} & \text{if} ~~ b_i \in A(X_i).
\end{cases}
$$
The morphism of simplicial coloured $\MM$-operads $(f,\phi)\colon \cO \to \cP$ can now be defined on colours by the inclusion $f \colon A \hookrightarrow A \sqcup F=B$ and by the identity $\phi \colon \cO \to i^*(\cP)=\cO$. In particular, the morphism $(f,\phi)$ is fully faithful. Note also that if $u\colon K \to H$ is an equivalence of categories, $(f,\phi)$ is moreover an equivalence. To conclude the proof, we must show that $(f,\phi)$ is the desired pushout of the morphism $i_!(\Id{F}\times u)$. This can be done in an analogous way as in the proof of \cite[Lemma 1.29]{cismo}.
\end{proof}
\end{proposition}

\begin{proposition}\label{Jequi}
The morphisms in $\overline{J}$ are weak equivalences, where $\overline{J}$ denotes the saturation of $J$.
\begin{proof}
Using the small object argument, we know that morphisms in $\overline{J}$ are retracts of sequential colimits of pushouts of morphisms in $J$. The class of weak equivalences in $\soperc{\MM}$ is closed under retracts, and, using that pushouts of morphisms in $J$ are monomorphisms on colours along with \ref{filcol}, we just have to prove that pushouts of morphisms in $J$ are weak equivalences. Consider first pushouts along morphisms $(\Id{A},\phi)\colon (A,\cO) \to (A,\cP)$ in $J_A$. In this case, the pushout of $(\Id{A},\phi)$ along a morphism $(g,\psi)\colon (A,\cO) \to (B,\cQ)$ is the morphism $(\Id{B},\eta)$ where $\eta$ is obtained via the following pushout 
$$
\begin{tikzcd}
g_! \cO \arrow[r,"\psi'"] \arrow[d,"g_1 \phi "'] & \cQ \arrow[d,"\eta"] \\
g_!\cP \arrow[r] & \cR
\end{tikzcd}
$$
in the category $\soper{\MM}^B$ where $\psi'\colon g_! \cO \to \cQ$ is obtained from $\psi$ by adjunction. Since $g_!$ is a left Quillen functor and because $\phi$ is by definition a trivial cofibration in $\soper{\MM}^A$, it follows that $\eta$ is a trivial cofibration in $\soper{\MM}^B$ and, in particular, $(\Id{B},\eta)$ is a local weak equivalence in $\soperc{\MM}$. The morphism $(\Id{B},\eta)$ is also trivially essentially surjective, and it follows that it is a weak equivalence in $\soperc{\MM}$.
Consider now a morphism $(f,\phi)\colon (A,\cO) \to (B,\cP)$ that belongs to $J_0$.  There exists $X$ in $\MM_0$ and a simplicial category $H$ satisfying the conditions of $(A_2)$ such that $(f,\phi)$ is the image under $i_!$ of the associated morphism $\hom(-,X) \to H\times \hom(-,X)$ in $\fun(\opo{\MM_0},\sCat)$. Consider now a pushout square
\begin{equation}\label{pushoutth}\tag{C}
\begin{tikzcd}
i_! (\hom(-,X)) \arrow[r] \arrow[d] & (C,\cQ) \arrow[d,"(g\text{,}\eta)"] \\
i_!(H \times \hom(-,X)) \arrow[r] & (D,\cR)
\end{tikzcd}
\end{equation}

in $\soperc{\MM}$ and let us show that $(g,\eta)$ is a weak equivalence in $\soperc{\MM}$.
As in the proof of \cite[Proposition A.25]{cismo}, the pushout square \ref{pushoutth} can be decomposed as follows 
$$
\begin{tikzcd}
i_! (\hom(-,X)) \arrow[r] \arrow[d]  & \arrow[dd,"(g\text{,}\eta)",bend left=40] (C,\cQ) \arrow[d,"(h\text{,}\nu)"'] \\
i_!(K \times \hom(-,X)) \arrow[d] \arrow[r] & (C,\cS) \arrow[d,"(k\text{,}\mu)"'] \\
i_!(H \times \hom(-,X)) \arrow[r] & (D,\cR)
\end{tikzcd}
$$
where $K$ is $u^*(H)$ with $u\colon \{0\} \hookrightarrow \{0,1\}$ and we have to show that both the morphisms $(h,\nu)\colon (C,\cQ) \to (C,\cS)$ and $(k,\mu)\colon (C,\cS) \to (D,\cR)$ are weak equivalences.
The morphism $i_!(\hom(-,X)) \to i_!(K\times \hom(-,X))$ is constant on objects and it follows that $(h,\nu)\colon (C,\cQ) \to (C,\cS)$ is a weak equivalence if $i_!(\hom(-,X)) \to i_!(K\times \hom(-,X))$ is a trivial cofibration in $\soper{\MM}^{\{0\}}$. Using \cite[Lemma 1.28]{cismo}, $K$ is cofibrant in $\sCat^{\{0\}}$ and it follows that $i_!(\hom(-,X)) \to i_!(K\times \hom(-,X))$ is a trivial cofibration. We now have to prove that the morphism $(k,\mu)\colon (C,\cS) \to (D,\cR)$ is a weak equivalence. Using \ref{simpob}, we can compute the pushout levelwise in the category $\fun(\opo{\Delta},\operc{\MM})$. Levelwise, the simplicial functor $K \to H$ satisfies the conditions of \ref{pushj} and we deduce that $(k,\mu)\colon (C,\cS) \to (D,\cR)$ is fully faithful (we even have isomorphisms of simplicial sets between the simplicial sets of operations). We now have to show that $(k,\mu)\colon (C,\cS) \to (D,\cR)$ is essentially surjective. The functor $\pi_0 \colon \soperc{\MM} \to \operc{\MM}$ preserves colimits and it follows that the square 
$$
\begin{tikzcd}
i_!(\pi_0(K)\times F) \arrow[r] \arrow[d] & \pi_0(D,\cS) \arrow[d,"\pi_0(k\text{,}\mu)"] \\
i_!(\pi_0(H)\times F) \arrow[r] &  \pi_0(D,\cR)
\end{tikzcd}
$$
is a pushout in $\operc{\MM}$, and, because $\pi_0(u)\colon \pi_0(K) \to \pi_0(H)$ is an equivalence, it follows from \ref{pushj} that $\pi_0 (k,\mu)$ is essentially surjective.
\end{proof}
\end{proposition}

Finally, we prove the main Theorem of this section regarding the existence of the model structure on simplicial coloured $\MM$-operads.
\begin{theorem}\label{mainTH}
The category of simplicial coloured $\MM$-operads $\soperc{\MM}$ is endowed with a cofibrantly generated model structure in which a morphism $(f,\phi)\colon (A,\cO) \to (B,\cP)$ is a weak equivalence if it is a fully faithful and essentially surjective morphism and in which a morphism is a fibration if is is an isofibration in the sense of \ref{weakfibop}.
\begin{proof}
We check the conditions of \cite[Theorem 2.1.19]{Hov} where the sets of generating cofibrations and trivial cofibrations are respectively $I$ and $J$. Note first that we have already proven that $\soperc{\MM}$ is bicomplete in \ref{bicomcol}. The first condition is proven using both \ref{23} and \ref{retract}. The conditions (2) and (3) were proven in \ref{small}. The conditions (5), (4), and (6) follow from \ref{isovsJ}, \ref{isowvsI}, and \ref{Jequi}. 
\end{proof}
\end{theorem}

\begin{example}\label{modparacop}
For every orbital pair $(\sO,\sT)$, we obtain a model structure on the category $\soperc{(\sO,\sT)}.$ In particular, if we consider the final category $*$, we find back the existence of the model structure on $\soperc{}$ as constructed in \cite[Theorem 1.14]{cismo} and \cite[Theorem 6]{Rob}.
\end{example}

\begin{example}
If we consider the $\FF$-multicategory $\und{\NN}$ of \ref{exnonsym}, we obtain a model structure on the category of simplicial multicategories.
\end{example}

\section{Complete \texorpdfstring{$\MM$}{M}-Segal spaces}\label{sec5}

In this final section, we define a notion of \emph{complete $\MM$-Segal spaces} that should produce the same homotopy theory as simplicial coloured $\MM$-operads. We define first, for every $\FF$-multicategory $\MM$, an algebraic pattern $\Delta_{\MM}^{\op}$ in the sense of \cite[Definition 2.1]{ChuHaug}, which will encode our $\MM$-Segal spaces.

\begin{definition}\label{defgendelta}
	Let $\DD$ be a double category. The category $\Del{\DD}$ is the (cocartesian) Grothendieck construction of the nerve $N(\DD)\colon \opo{\Delta} \to \Cat$ of $\DD$. 
	If $\MM$ is an $\FF$-multicategory, then we define $\Del{\MM}$ as the category $\Del{\FF \MM}$. 
    We have a cocartesian fibration $\pi_{\DD}\colon \opo{\Del{\DD}} \to \opo{\Delta}$ and it follows that the algebraic pattern structure of $\Delta^{\op}$, denoted by $\Delta^{\op,\natural}$ in \cite[Example 3.3]{ChuHaug}, induces the structure of an algebraic pattern on $\opo{\Delta_{\DD}}$. Explicitly, the structure of an algebraic pattern on $\opo{\Delta_{\DD}}$ is defined as follows:
\begin{enumerate}
	\item a morphism in $\Delta_{\DD}^{\op}$ is \emph{inert} if the underlying morphism in $\Delta^{\op}$ is inert;
	\item a morphism in $\Delta_{\DD}^{\op}$ is \emph{active} if it is $\pi_{\DD}$-cocartesian and if the underlying morphism in $\Delta^{\op}$ is active;
	\item an object $([n],f)$ in $\Delta_{\DD}^{\op}$ is elementary if $n=1$ or $n=0$.
\end{enumerate}
For every $\FF$-multicategory $\MM$, the category $\Delta_{\MM}^{\op}$ admits the structure of an algebraic pattern whose inert and active morphisms are the ones of $\Delta_{\FF \MM}^{\op}$ and where the elementary objects are the objects $([1],f)$ where $f$ belongs to $\MM_1$ and $([0],x)$ where $x$ belongs to $\MM_0$.
\end{definition}

\begin{remark}
Note that, in \ref{defgrotcons}, we were considering the cartesian Grothendieck construction, whereas in the previous definition, we consider the cocartesian one.
\end{remark}

\begin{definition}
An \emph{$\MM$-Segal space} is a $\Delta_{\MM}^{\op}$-Segal object in $\sS$, in the sense of \cite[Definition 2.7]{ChuHaug}. We denote by $\Seg{\MM}(\sS)$ the $\infty$-category of $\MM$-Segal spaces. 
\end{definition}

We define a notion of completeness for $\MM$-Segal spaces. Every object $X$ of $\MM$ provides a section $i_X\colon \Delta \to \Delta_{\MM}$ of $\pi_{\MM}$. Moreover, the functor $i_X$ is a Segal morphism in the sense of \cite[Definition 4.2]{ChuHaug} and it induces a functor $i^*_X\colon \Seg{\MM}(\sS) \to \Seg{}(\sS)$.

\begin{definition}
An $\MM$-Segal space $\cO$ is \emph{complete} if for every $X$ in $\MM_0$, the Segal space $i_X^*(\cO)$ is complete. We denote by $\CSeg{\MM}(\sS)$ the full subcategory of $\Seg{\MM}(\sS)$ spanned by complete $\MM$-Segal spaces.
\end{definition}

\begin{example}
For the $\FF$-multicategory $\und{\Sigma_*}$ of \ref{exsymmop}, we obtain a slight modification of the algebraic pattern $\Delta_{\Fin}$ defined by Barwick in \cite{Bar2}. More precisely, the category $\Delta_{\Fin}$ is a wide subcategory of our $\Delta_{\und{\Sigma_*}}$ and this inclusion satisfies the conditions of \cite[Corollary 2.64]{Bark1}. This implies that our $\und{\Sigma_*}$-Segal spaces correspond to the $\Fin$-Segal spaces of Barwick.
\end{example}

\begin{example}
For every orbital pair $(\sO,\sT)$, we obtain a notion of \emph{complete $(\sO,\sT)$-Segal spaces} using the $\FF$-multicategory $\und{\Sigma_{(\sO,\sT)}}$ of \ref{defparaop}. We denote by $\CSeg{(\sO,\cT)}$ the associated $\infty$-category.
\end{example}

In the same way simplicial categories can be compared with complete Segal spaces, we make the following conjecture.

\begin{conjecture}\label{segvssimp}
	The $\infty$-category obtained from $\soperc{\MM}$ by inverting fully faithful essentially surjective morphisms (defined in \ref{weakfibop}) is equivalent to the $\infty$-category $\CSeg{\MM}(\sS)$ of complete $\MM$-Segal spaces. 
\end{conjecture}

\begin{definition}
For every orbital pair $(\sO,\sT)$, we denote by $\infty$-$\oper{(\sO,\sT)}$ the $\infty$-category of fibrous patterns over $\Span_{\sT}\left(\FF \sO\right)$ in the sense of \cite{Nat2}. The objects of $\infty$-$\oper{(\sO,\sT)}$ will be refered to as \emph{$(\sO,\sT)$-$\infty$-operads}. Note that if $\sO$ is an atomic orbital category, then it follows from \cite[Theorem A.1]{Nat2} that the $\infty$-category $\infty$-$\oper{(\sO,\sO)}$ is equivalent to the $\infty$-category of $\sO$-$\infty$-operads of Nardin and Shah \cite{DJ}.
\end{definition}

A direct consequence of \ref{segvssimp} would be that the $\infty$-category associated with $\soperc{(\sO,\sT)}$ can be compared with complete $(\sO,\sT)$-Segal spaces.

\bibliographystyle{amsalpha}
\bibliography{references}

\end{document}